# The Ordered Distribution of Natural Numbers on the Square Root Spiral

- Harry K. Hahn -

Germany

mathematical analysis by

- Kay Schoenberger -

Humboldt-University
Berlin

20. June 2007

## Abstract :

Natural numbers divisible by the same prime factor lie on defined spiral graphs which are running through the "Square Root Spiral" ( also named as "Spiral of Theodorus" or "Wurzel Spirale" or "Einstein Spiral" ). Prime Numbers also clearly accumulate on such spiral graphs.
And the square numbers 4, 9, 16, 25, 36 … form a highly three-symmetrical system of three spiral graphs, which divide the square-root-spiral into three equal areas.
A mathematical analysis shows that these spiral graphs are defined by quadratic polynomials.
The Square Root Spiral is a geometrical structure which is based on the three basic constants: 1, sqrt2 and π (pi) , and the continuous application of the Pythagorean Theorem of the right angled triangle.
Fibonacci number sequences also play a part in the structure of the Square Root Spiral. Fibonacci Numbers divide the Square Root Spiral into areas and angle sectors with constant proportions. These proportions are linked to the "golden mean" ( golden section ), which behaves as a self-avoiding-walk-constant in the lattice-like structure of the square root spiral.

## Contents of the general section — Page

## Contents of the mathematical section

# General Section

## 1 Introduction to the Square Root Spiral :

The Square Root Spiral ( or "Wheel of Theodorus" or "Einstein Spiral" or "Wurzel Spirale" ) is a very interesting geometrical structure, in which the square roots of all natural numbers have a clear defined spatial orientation to each other. This enables the attentive viewer to find many interdependencies between natural numbers, by applying graphical analysis techniques. Therefore, the square root spiral should be an important research object for all professionals working in the field of number theory !

Here is a first impressive image of the Square Root Spiral :

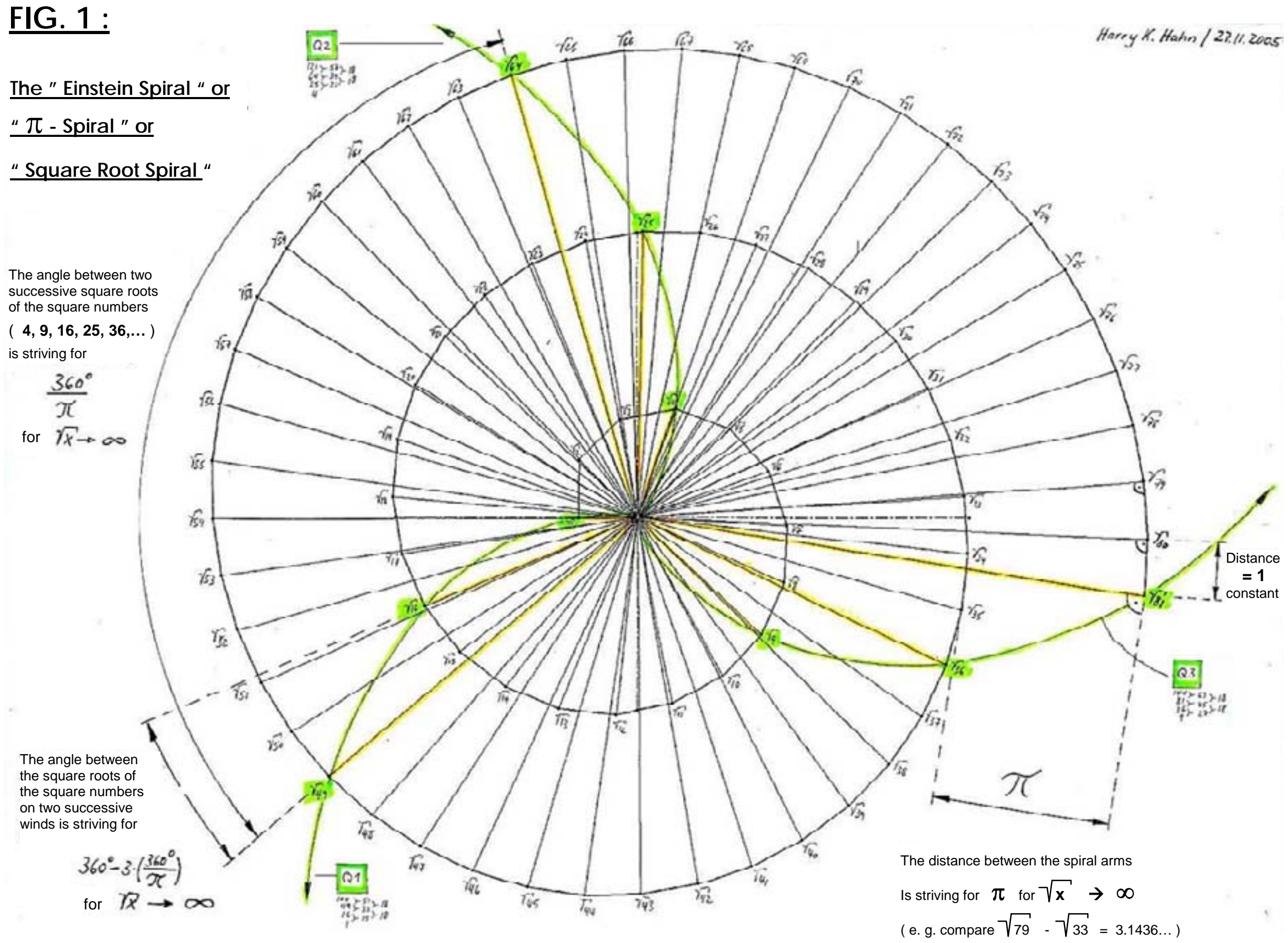

The most amazing property of the square root spiral is surely the fact that the distance between two successive winds of the Square Root Spiral quickly strives for the well known geometrical constant **π** !!

Mathematical proof that this statement is correct is shown in **Chapter 1 " The correlation with π "** in the mathematical section.

→ **Table 1** in the **Appendix** shows an approximate analysis of the development of the distance between two successive winds of the Square Root Spiral ( or Einstein Spiral ).
In principle this analysis uses the length difference of two " square root rays" which differ by nearly exactly one wind of the square root spiral to each other.
( see example on **FIG. 1 :** sqrt 79 – sqrt 33 = **3.1436…** )

Another striking property of the Square Root Spiral is the fact, that the square roots of all square numbers ( 4, 9, 16, 25, 36… ) lie on 3 highly symmetrical spiral graphs which divide the square root spiral into three equal areas ( see FIG.1 : graphs **Q1, Q2** and **Q3** drawn in green ). For these three graphs the following rules apply :

1.) The angle between successive Square Numbers ( on the "Einstein Spiral" ) is striving for 360 °/π for sqrt( X ) → ∞

2.) The angle between the Square Numbers on two successive winds of the "Einstein Spiral" is striving for 360 ° - 3x( 360°/π ) for sqrt( X ) → ∞

Proof that these propositions are correct shows **Chapter 2 " The Spiral Arms"** in the mathematical section.

**The Square Root Spiral** develops from a right angled base triangle ( **P1** ) with the two legs ( cathets ) having the length 1, and with the long side ( hypotenuse ) having a length which is equal to the square root of 2.

→ **see FIG. 2 and 4**

The square root spiral is formed by further adding right angled triangles to the base triangle **P1** ( see FIG 4) In this process the longer legs of the next triangles always attach to the hypotenuses of the previous triangles. And the longer leg of the next triangle always has the same length as the hypotenuse of the previous triangle, and the shorter leg always has the length 1.
In this way a spiral structure is developing in which the spiral is created by the shorter legs of the triangles which have the constant length of 1 and where the lengths of the radial rays ( or spokes ) coming from the centre of this spiral are the square roots of the natural numbers ( sqrt 2 , sqrt 3, sqrt 4, sqrt 5 …. ).

→ **see FIG. 1 and 4**

The special property of this infinite chain of triangles is the fact that all triangles are also linked through the Pythagorean Theorem of the right angled triangle. This means that there is also a logical relationship between the imaginary square areas which can be linked up with the cathets and hypotenuses of this infinite chain of triangles ( → all square areas are multiples of the base area 1 , and these square areas represent the natural numbers N = 1, 2, 3, 4,…..) → **see FIG. 2 and 3**. This is an important property of the Square Root Spiral, which might turn out someday to be a "golden key" to number theory !

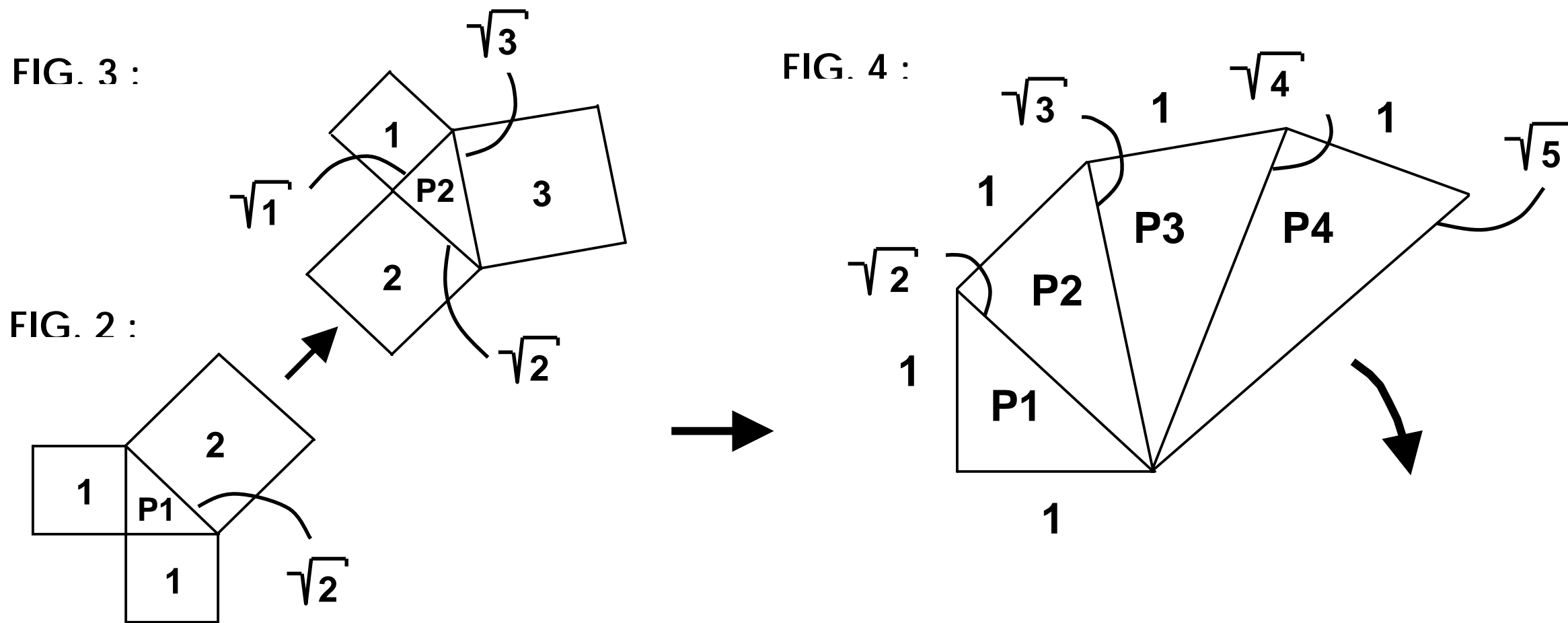

FIG. 1 shows the further development of the square root spiral or "Einstein Spiral" if one rectangular triangle after the other is added to the growing chain of triangles as described in FIG. 4.

For my further analysis I have created a square root spiral consisting of nearly 300 precise constructed triangles. For this I used the CAD Software SolidWorks. The length of the hypotenuses of these triangles which represent the square roots from the natural numbers 1 to nearly 300, has an accuracy of 8 places after the decimal point. Therefore, the precision of the square root spiral used for the further analysis can be considered to be very high. ( a bare Square Root Spiral can be found in the Appendix → see FIG. 17 )

The lengths of the radial rays ( or spokes ) coming from the centre of the square root spiral represent the square roots of the natural numbers ( n = { 1, 2, 3, 4,...} ) in reference to the length 1 of the cathets of the base triangle P1 ( see FIG. 4 ). And the natural numbers themselves are imaginable by the areas of "imaginary squares", which stay vertically on these "square root rays". → **see FIG. 5** (compare with FIG.3 )

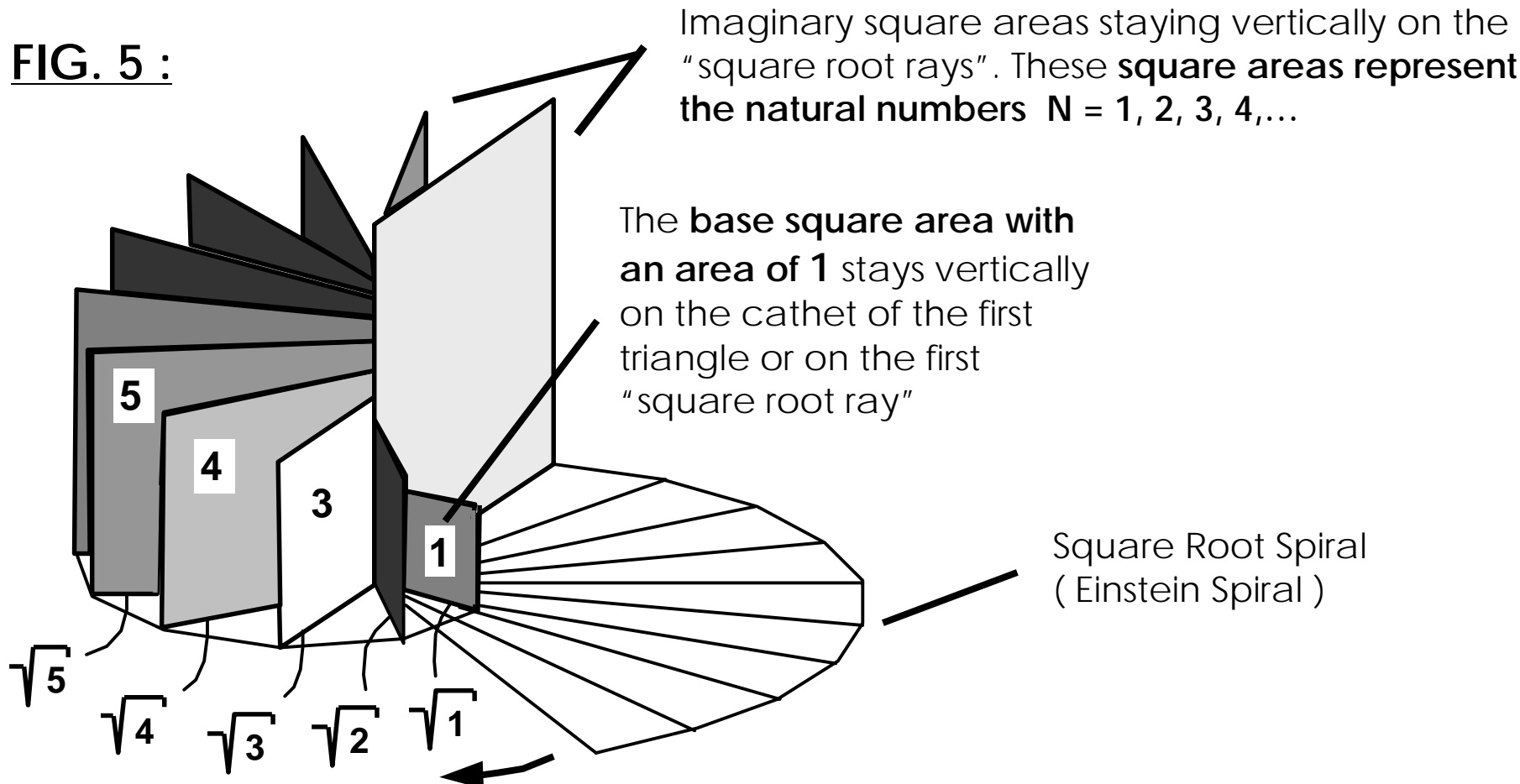

→ The "square root rays" of the Einstein Spiral can simply be seen as a projection of these spatially arranged "imaginary square areas", shown in FIG. 5, onto a two-dimensional plane.

## 2 Mathematical description of the Square Root Spiral

Comparing the Square Root Spiral with different types of spirals ( e.g. logarithmic-, hyperbolic-, parabolic- and Archimedes- Spirals ), then the Square Root Spiral obviously seems to belong to the Archimedes Spirals.

An Archimedes Spiral is the curve ( or graph ) of a point which moves with a constant angle velocity around the centre of the coordinate system and at the same time with a constant radial velocity away from the centre. Or in other words, the radius of this spiral grows proportional to its rotary angle.

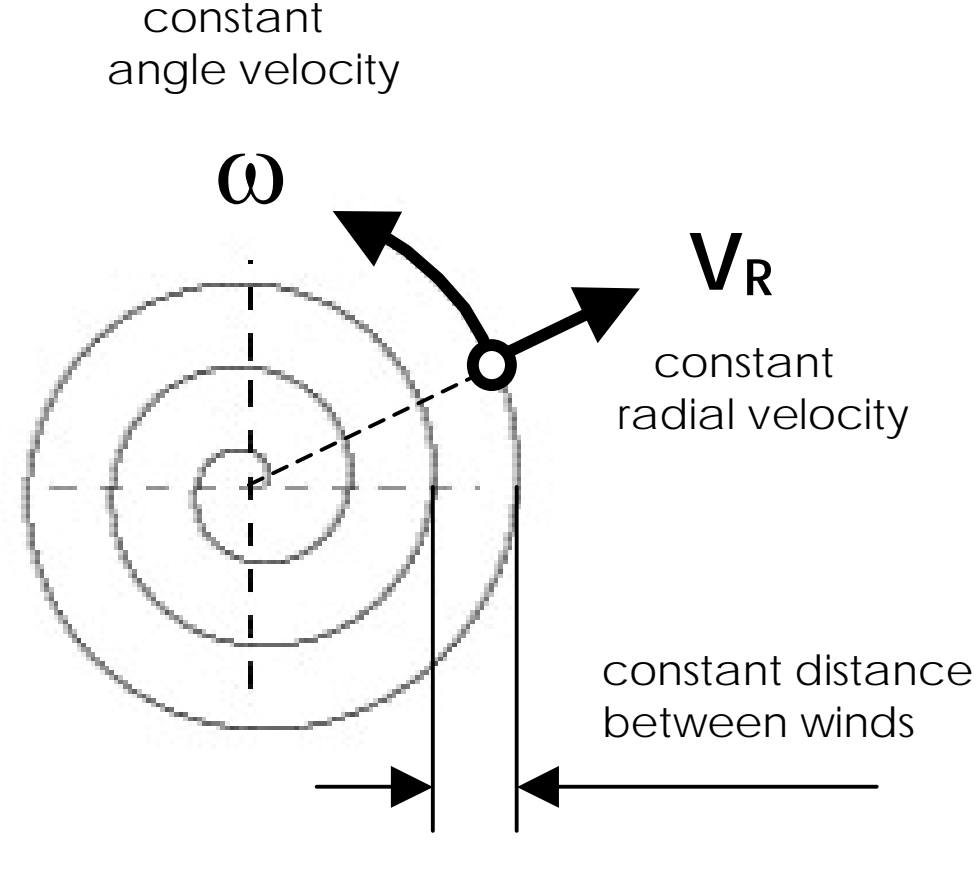

**Archimedes Spiral**

In polar coordinate style the definition of an Archimedes Spiral reads as follows :

$$r(\varphi) = a\varphi \qquad \text{with} \qquad a = \text{const.} = \frac{V_R}{\omega} > 0$$

for $r \rightarrow \infty$ the Square Root Spiral is an Archimedes Spiral with the following definition :

$$r(\varphi) = a\varphi + b$$

with $a$ = const. and $b$ = const.

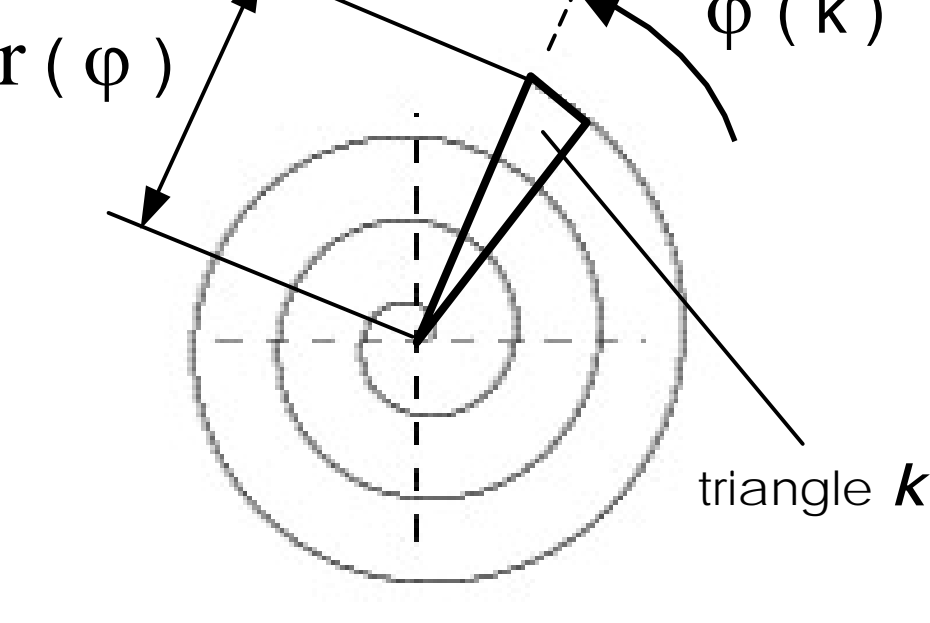

The values of the parameters $a$ and $b$ are

$$a = \frac{1}{2} \quad \text{and} \quad b = -\frac{c_2}{2} \quad ; \text{ with}$$

$c_2$ = Square Root Spiral Constant

$c_2$ = - 2.157782996659….

Hence the following formula applies for the Square Root Spiral :

$$r(\varphi) = \frac{1}{2}\varphi + 1.078891498\ldots.. \qquad \text{for } r \rightarrow \infty$$

for $r \rightarrow \infty$ therefore the growth of the radius of the Square Root Spiral after a full rotation is striving for $\pi$ ( corresponding to the angle of a full rotation which is $2\pi$ )

Note : The mathematical definitions shown on this page and on the following page can also be found either in the mathematical section of this paper, or in other studies referring to the Square Root Spiral. → e.g. a mathematical analysis of the Square Root Spiral is available on the following website : http://kociemba.org/themen/spirale/spirale.htm

→ Also note that in the mathematical section of my paper ( contributed by Mr Kay Schoenberger ) $t_n$ is used instead of $\varphi_n$ and $\omega(k)$ instead of $\varphi(k)$

**Further dependencies in the Square Root Spiral :**

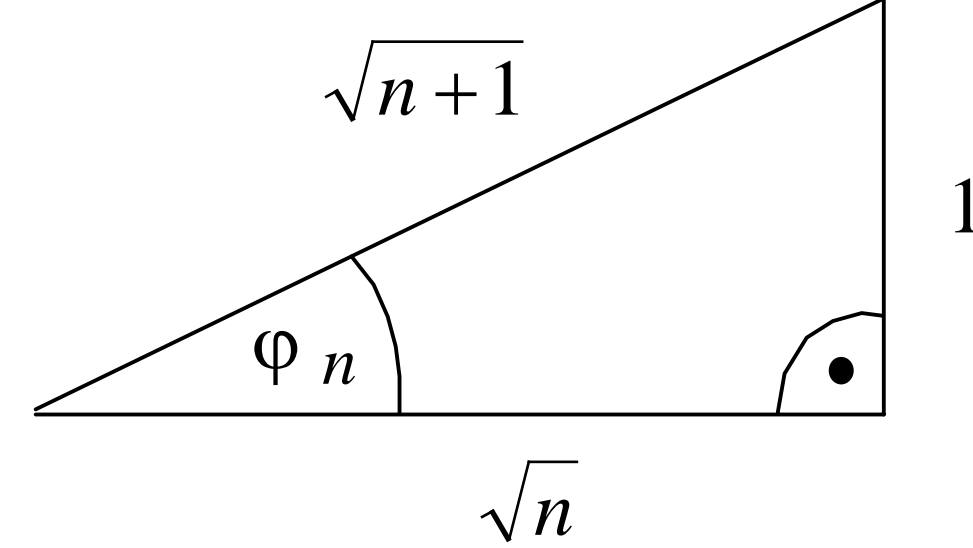

If $\varphi_n$ is the angle of the ***n***th spiral segment ( or triangle ) of the Square Root Spiral, then

$\tan(\varphi_n) = \frac{1}{\sqrt{n}}$ ; ( ratio $\frac{\text{counter cathet}}{\text{cathet}}$ )

If the ***n***th triangle is added to the Square Root Spiral the growth of the angle is

$$\varphi_n = \arctan\left(\frac{1}{\sqrt{n}}\right)$$

; <u>Note</u> : angle in radian

The total angle $\varphi(k)$ of a number of **k** triangles is

$$\varphi(k) = \sum_{n=1}^{k} \varphi_n$$

or described by an integral

$$\int_0^k \arctan\left(\frac{1}{\sqrt{n}}\right) dn + c_1(k)$$

$$\Rightarrow \quad \varphi(k) = 2\sqrt{k} + c_2(k) \quad \text{with} \quad \lim_{k\to\infty} c_2(k) = \text{const.} = -2.157782996659\ldots..$$

$c_2$ = Square Root Spiral Constant

The growth of the radius of the Square Root Spiral at a certain triangle $n$ is

$$\Delta r = \sqrt{n+1} - \sqrt{n}$$

The radius $r$ of the Square Root Spiral ( i.e. the big cathet of triangle **k** ) is

$r(k) = \sqrt{k}$ and by converting the above shown equation for $\varphi(k)$ it applies that

$$r(k(\varphi)) = r(\varphi) = \sqrt{\frac{1}{4}(\varphi - c_2(\varphi))^2} = \frac{1}{2}\varphi - \frac{c_2}{2}$$

For large $n$ it also applies that $\varphi_n$ is approximately $\frac{1}{\sqrt{n}}$ and $\Delta r$ has pretty well half of this value, that is $\frac{1}{2\sqrt{n}}$ , what can be proven with the help of a Taylor Sequence.

## 3 The distribution of the square numbers 4, 9, 16, 25, 36, ... on the Square Root Spiral :

The square roots of the square numbers (1), 4, 9, 16, 25, 36, 49,… lie in 3 areas which are arranged highly symmetrically around the center of the square root spiral.
Here the square numbers themselves can be represented by the mentioned imaginary square areas which stay vertically on the "square root rays"
→ **see FIG. 6**

And the square roots of the square numbers, which are the numbers 1, 2, 3, 4, 5, 6,… are the "square root rays" which form the base lines of these imaginary quadratic areas .

Only the square roots of the square numbers are whole numbers or natural numbers.
That's why the 3-symmetrical distribution of these numbers on the square root spiral must have an important meaning !

Especially if we consider that the Square Root Spiral is precisely divided in 3 equal areas by the square numbers !

But it seems nobody has yet taken notice of this amazing fact and tried to explain it !!

**FIG. 6 :**

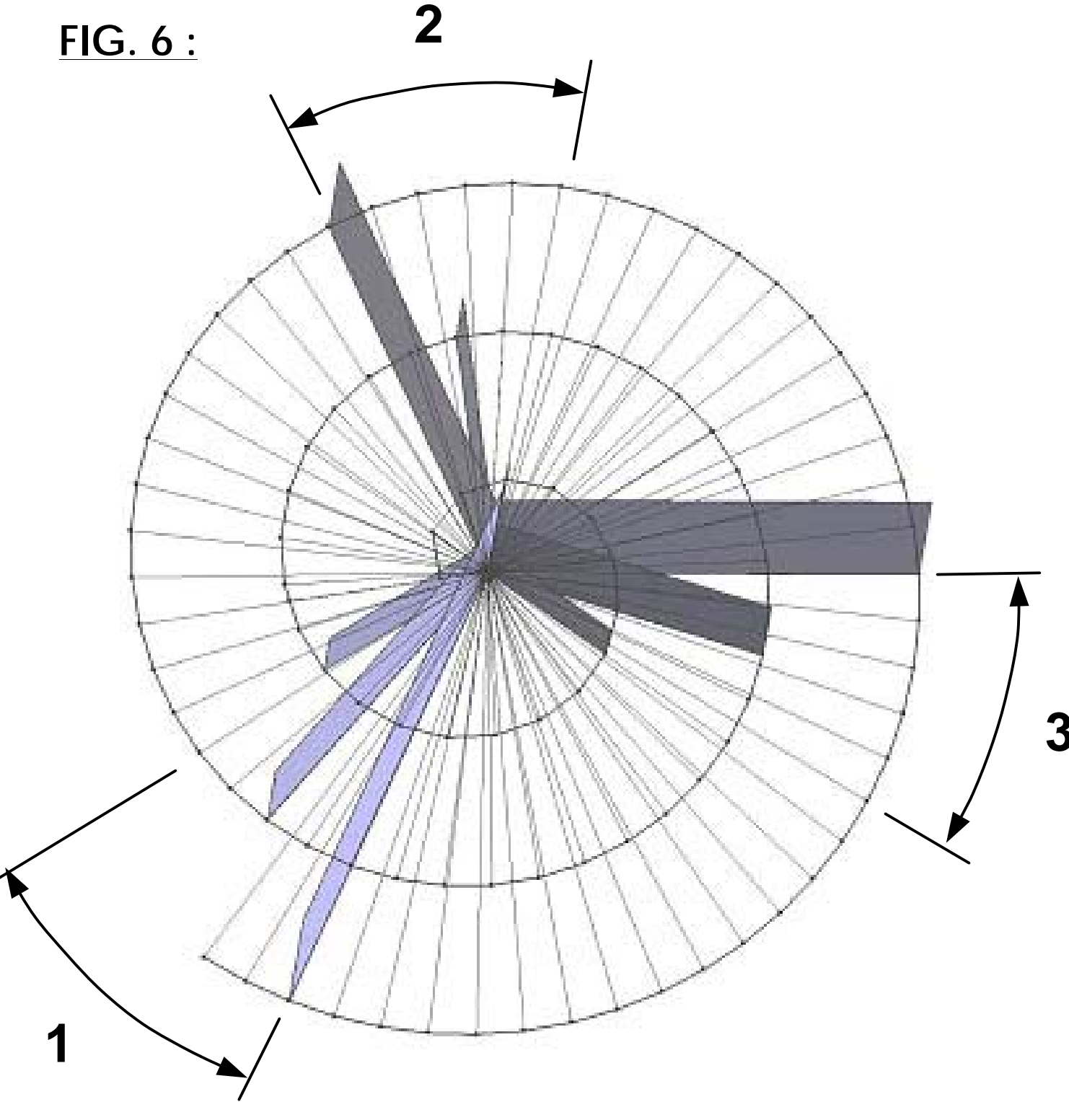

The " square root rays" of the square numbers are arranged in a way that their outer ends lie on three spiral graphs ( quadratic polynomials ) as shown in FIG. 1 ( → spiral graphs drawn in green ).

### 3.1 Listing of important properties of the three spiral graphs containing the square numbers :

- **The Square Numbers** lie on **3** highly symmetrical spiral graphs with a **positive** rotation direction (drawn in green ). → **see FIG.1** → These 3 spiral-graphs are defined by the following :
  **3 Quadratic Polynomials :** **Q1** = $9x^2 + 6x + 1$
  **Q2** = $9x^2 + 12x + 4$
  **Q3** = $9x^2 + 18x + 9$

  ( → see also Table 3-B at page 34 in the Appendix ! )

- The **3** spiral graphs **Q1 – Q3** are arranged in an angle of around 120° to each other ( referring to the center of the Square Root Spiral )

- It applies : **Q1** contains the square number sequence **1, 16, 49, 100, 169,…**
  ( the square roots of these numbers are : 1, 4, 7, 10, 13…→ difference = 3 )

  **Q2** contains the square number sequence **4, 25, 64, 121, 196,…**
  ( the square roots of these numbers are : 2, 5, 8, 11, 14,…→ difference = 3 )

  **Q3** contains the square number sequence **9, 36, 81, 144, 225,…**
  ( the square roots of these numbers are : 3, 6, 9, 12, 15,…→ difference = 3 )

  ( → in the **Q3** - sequence **all numbers are also divisible by 3** ! )

→ **FIG. 15** at page 30 in the **Appendix** shows the exact geometry of spiral graph **Q1** :

- The angle between successive square numbers on the Square Root Spiral ( "Einstein Spiral" ) is striving for 360 °/π for sqrt( X ) → ∞ → see **FIG. 1 / FIG. 8 & mathematical section** )

- The angle between the square numbers on two successive winds of the Square Root Spiral ( "Einstein Spiral" ) is striving for 360 ° - 3x( 360°/π ) for sqrt( X ) → ∞ → see **FIG. 1 / FIG. 8**

- Calculating the differences of the consecutive square numbers lying on one of the three spiral arms, and then further calculating the differences of these differences ( → " **2. Differential** " ), results in the constant value **18** for the three spiral graphs ( quadratic polynomials ) **Q1 – Q3**. ( **→ see difference values in FIG.1** beside the names of the spiral-graphs **Q1 – Q3** )

- The **3** spiral graphs containing the square numbers divide the square root spiral exactly into **3 equal ares.**

  **Proof** that this proposition is correct can be found in **Chapter 3 "Area equality"** in the **mathematical section**.

  The following analysis can also be used as a first approximate proof that this proposition is correct :

First we calculate the areas which lie between the square roots of the square numbers.

→ see the first three such areas on the square root spiral marked in green, yellow and red in FIG. 7

Then we always calculate the ratio of two such successive areas.

→ see calculation process shown below the image :

For $\sqrt{x}$ → ∞ the resulting ratio is striving for the value of **1** at infinity.

This first approximation indicates that the square root spiral is precisely divided into **3** equal areas by the square numbers !

**FIG. 7 :**

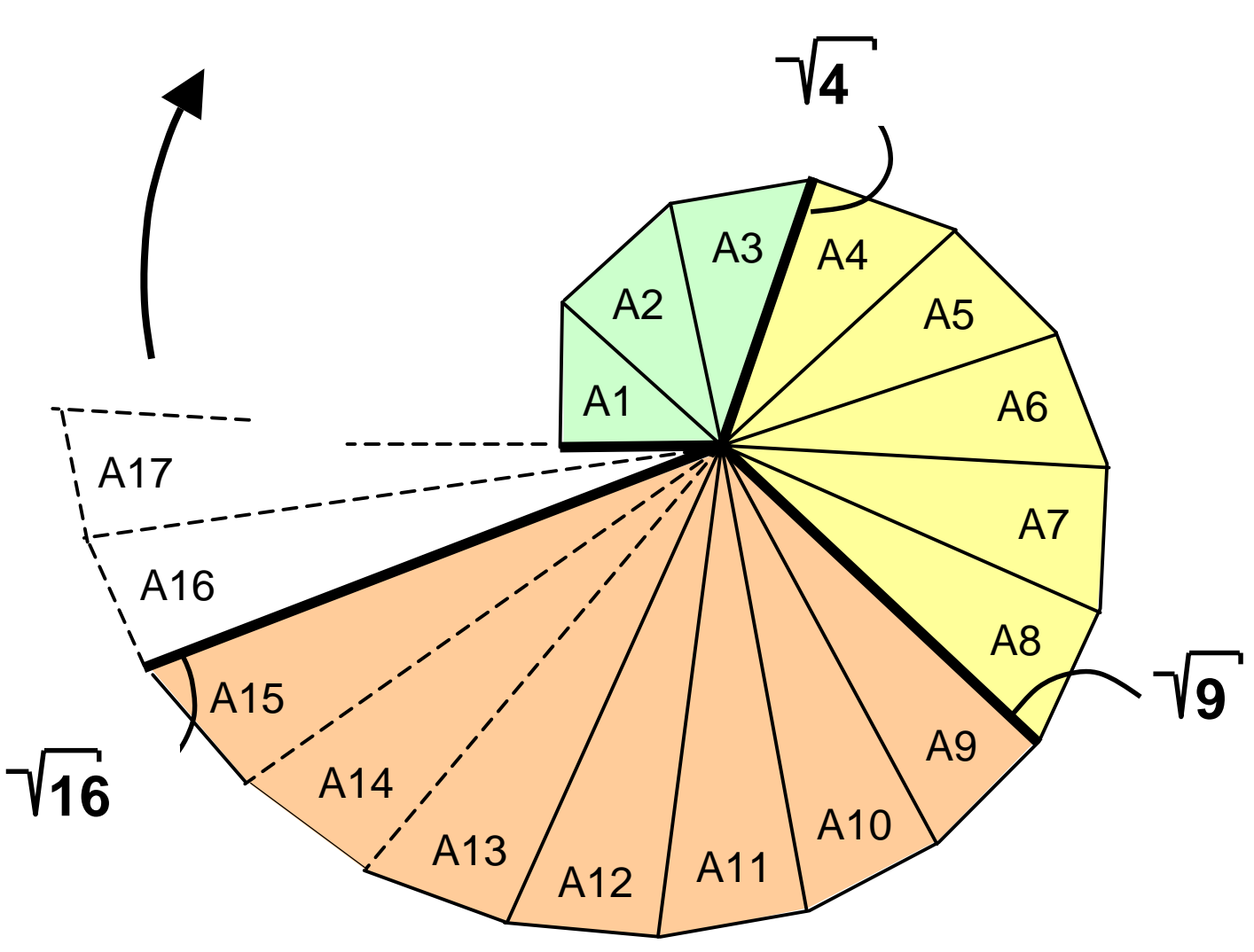

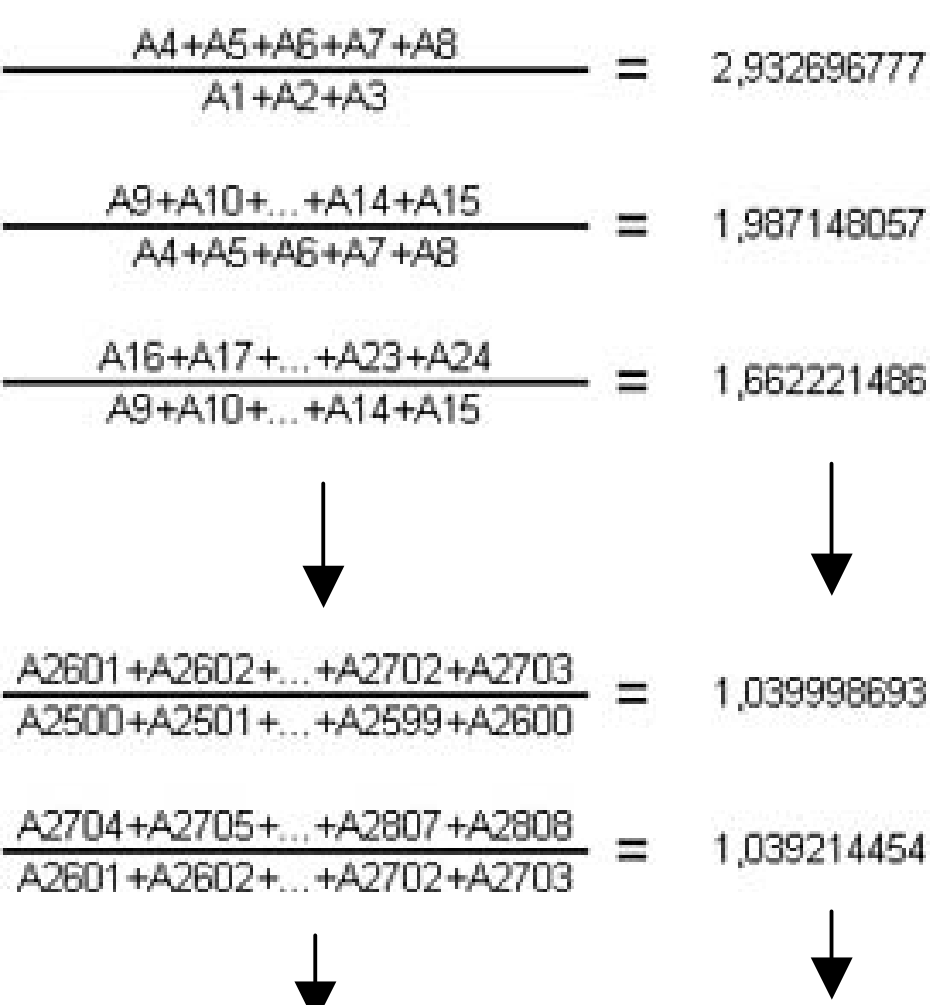

For the angles between the " square root rays" of the square numbers a similar approximation can be made as for the areas . → **see FIG. 1** and **FIG. 8 :**

**FIG. 8** shows the development of the angles between the "square root rays" of the square numbers.

The precisely measured angles indicate the correctness of the two following statements :

1.) The angle between successive square numbers on the Square Root Spiral ( "Einstein Spiral" ) is striving for 360 °/π for sqrt( X ) → ∞

2.) The angle between the square numbers on two successive winds of the Square Root Spiral ( "Einstein Spiral" ) is striving for 360 ° - 3x ( 360°/π ) for sqrt( X ) → ∞

**Proof** that these propositions are correct is shown in **Chapter 2 "The Spiral Arms"** (→ **mathematical section**)

For further mathematical analysis of the spiral-graphs Q1, Q2 and Q3 shown in FIG. 1 , I included the exact geometry of the spiral graph **Q1**, which contains the square roots of the square numbers 1, 16, 49, 100, 169,…

This graph together with accurate polar coordinates can be found in the **Appendix** → **see FIG. 15**

## FIG. 8 :

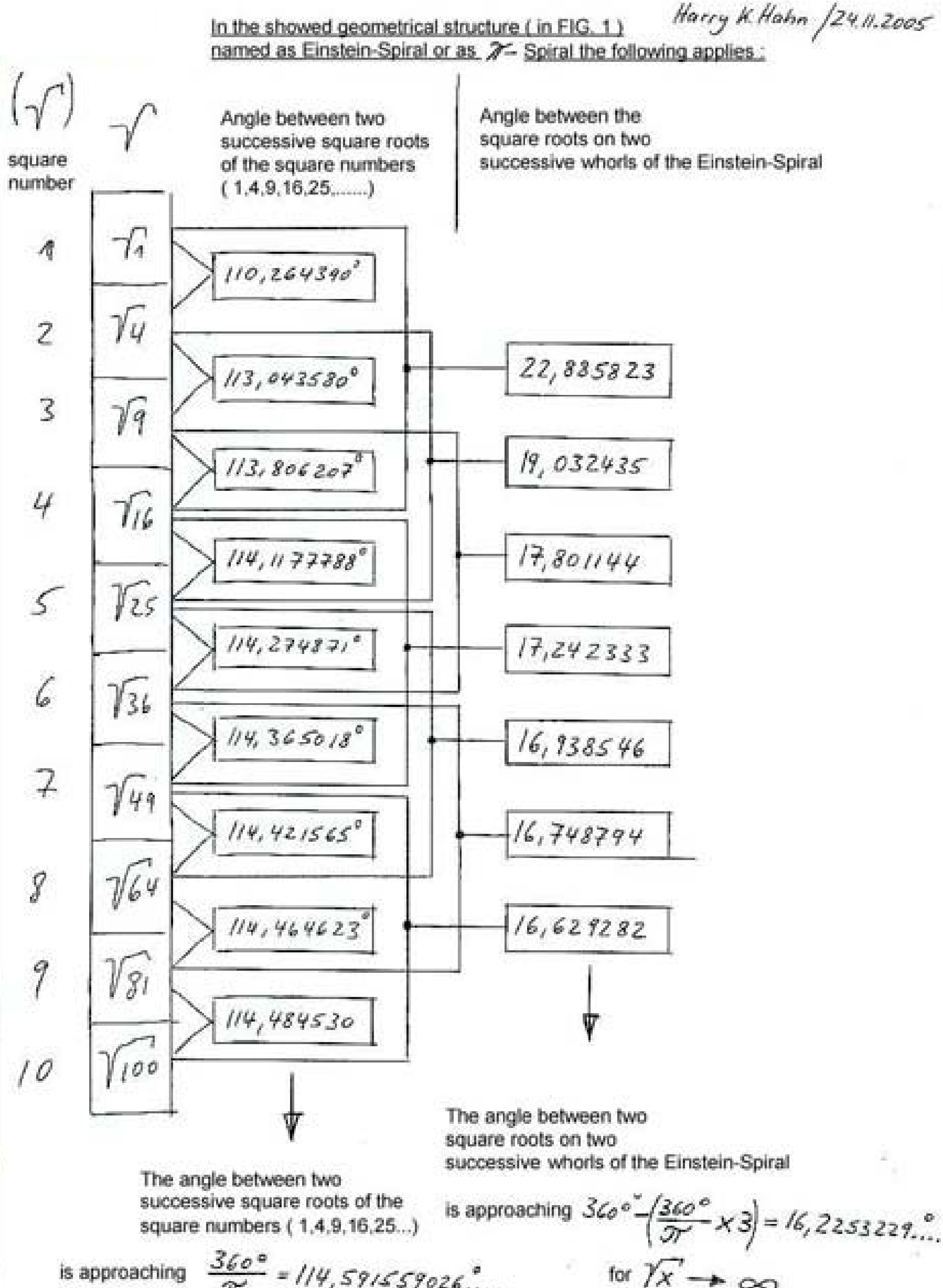

## 4 The distribution of natural numbers divisible by the prime factors 2, 3, 5, 7, 11, 13, 17,…

In comparison to the square numbers, which lie on three single spiral arms, which are symmetrically arranged around the center of the Square Root Spiral, all other natural numbers lie on "spiral graph systems" which consist of more than one spiral arm.
Here the natural numbers divisible by the prime factors 2, 3, 5, 7, 11 lie on more than one of these mentioned "spiral graph systems" with either a positive rotation direction or a negative rotation direction respectively. Natural numbers divisible by the prime factor 13 lie on only one spiral graph system with a positive rotation direction, but on two spiral graph systems with a negative rotation direction. And all natural numbers divisible by prime factors $\geq 17$ lie on only one spiral graph system with either a positive- or a negative rotation direction.

The following image **FIG. 9** shows for example the distribution of the natural numbers divisible by **11** on the Square Root Spiral. Here all numbers divisible by 11 are marked in yellow.

### 4.1 The distribution of natural numbers divisible by the prime factor 11 :

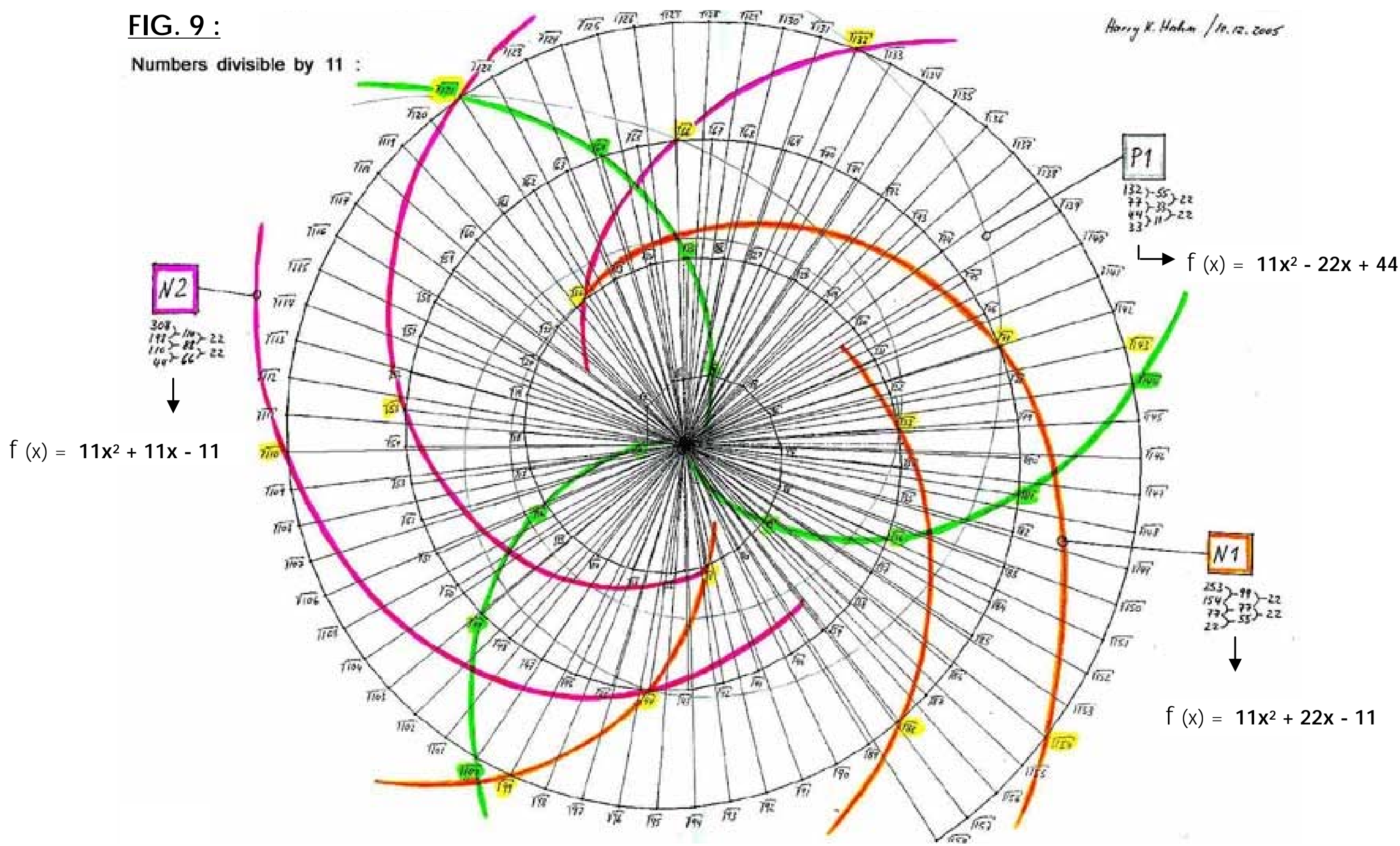

As mentioned before,…If I am now talking about the arrangement of the numbers divisible by 11 on the Square Root Spiral, I am actually referring to imaginary square areas, which stay vertical on certain radial rays of the square root spiral. The natural numbers divisible by 11, are represented by these imaginary square areas ( as explained in chapter 1 ). However in this analysis we only consider the projections of these imaginary square areas ( = radial rays ) onto a two-dimensional plane, for simplification.

From the image **FIG.9** it is evident that the ( square roots of the ) natural numbers divisible by 11 ( marked in yellow ) lie on defined spiral graphs which have their starting point in or near the centre of the Square Root Spiral. These spiral graphs have either a positive or a negative rotation direction.
A spiral graph which has a clockwise rotation direction shall be called negative (N) and a spiral graph which has a counterclockwise rotation direction shall be called positive (P).
The green spiral graphs show the three spiral-graphs which contain the square numbers 4, 9, 16, 25, 36, … which are drawn for reference only !

## 4.11 Properties of the spiral graph systems containing the numbers divisible by 11 :

- **The numbers divisible by 11** lie on **2** spiral-graph-systems with a **negative** rotation direction ( drawn in orange and pink ) and on **2** spiral-graph-systems with a **positive** rotation direction ( only one system drawn in light grey lines ! ). The **2** negative spiral graph systems are named **N1** and **N2** and the **2** positive systems are named **P1** and **P2** ( only P1 is drawn ! ).

  → Note : Not all spiral arms of the described spiral-graph-systems are drawn !

- The spiral-graph-systems **N1** and **N2** as well as **P1** and **P2** lie approximately point-symmetrical to each other ( in reference to the centre of the square root spiral )

- Calculating the differences of the consecutive numbers lying on one of the drawn spiral arms, and then further calculating the differences of these differences ( **2. Differential** ), results in the constant value **22** for the positive as well as the negative rotating spiral-graphs.

  → see difference calculation for **3 exemplary spiral arms** ( one spiral arm per system ) in FIG. 9 beside the names of the spiral-graph-systems N1, N2 and P1. ( → P2-system not shown ! )

  These **3** exemplary spiral arms are defined by the following quadratic polynomials :

  **Quadratic Polynomials of exemplary spiral arms :**

  | | | | |
  |---|---|---|---|
  | $\mathbf{11x^2 + 22x - 11}$ | $= 11\,(x^2 + 2x - 1)$ | → | belongs to **N1** – system |
  | $\mathbf{11x^2 + 11x - 11}$ | $= 11\,(x^2 + x - 1)$ | → | belongs to **N2** – system |
  | $\mathbf{11x^2 - 22x + 44}$ | $= 11\,(x^2 - 2x + 4)$ | → | belongs to **P1** - system |

  → **The following example shows how to calculate these quadratic polynomials.**

## 4.12 Calculation of the quadratic polynomial belonging to one exemplary spiral arm of the N1-system :

| | | | | | |
|---|---|---|---|---|---|
| Number sequence $N_1$ : | 22, | 77, | 154, | 253, … | ( → see number sequence beside the name of the **N1**-system ) |
| first difference : | 55 | 77 | 99 | | |
| second difference : | 22 | 22 | | | |
| third difference : | 0 | | | | |

Because the third differences are zero ( and this yields a quadratic polynomial ), we can use the short notation of the Newton Interpolation Polynomial to calculate the quadratic polynomial :

Here the following assignment is used :

| | | | | | |
|---|---|---|---|---|---|
| Number sequence: | $f(1)$ | $f(2)$ | $f(3)$ | $f(4)$ | … |
| first difference: | $f[1,2]$ | $f[2,3]$ | $f[3,4]$ | … | |
| second difference: | $2f[1,2,3]$ | $2f[2,3,4]$ | … | | |

with the short notation of the Newton Interpolation Polynomial we have the polynomial :

$$N(t) = f_1 + (t-1)f[1,2] + (t-1)(t-2)f[1,2,3]$$
$$= f_1 + (t-1)(f_2 - f_1) + \tfrac{1}{2}(t-1)(t-2)(f_1 - 2f_2 + f_3)$$

The generator polynomial for $N_1$ is therefore :

$$N_1(t) = 22 + (t-1)(77-22) + \tfrac{1}{2}(t-1)(t-2)(22 - 2\cdot 77 + 154)$$
$$= 11(t^2 + 2t - 1)$$

or in the general form of quadratic polynomials :

$$f(x) = 11\,(x^2 + 2x - 1) = \mathbf{11x^2 + 22x - 11}$$

→ Referring to the general quadratic polynomial $f(x) = \mathbf{ax^2 + bx + c}$ the following rules apply for the quadratic polynomials, belonging to the shown spiral-graphs :

Rules for coefficients **a, b** and **c** : ( or sequence of coefficients )

- **a** → equivalent to the **"2.Differential"** divided by **2**
- **b** → this coefficient ( or sequence of coefficients ) indicates the system of spiral-graphs it belongs to.
- **c** → describes the consecutive parallel distance of the spiral-graphs in the same system

→ **Please refer to chapter 2 "The Spiral Arms" in the mathematical section for a detailed mathematical explanation of the spiral arms ( or spiral-graphs ) shown in FIG. 1 / 9 / 10 / 11 / 12**

### 4.2 The distribution of natural numbers divisible by the prime factor 7 :

As a further example the following image **FIG.10** shows the distribution of the natural numbers divisible by **7** on the square root spiral. Here all numbers divisible by **7** are marked in yellow.

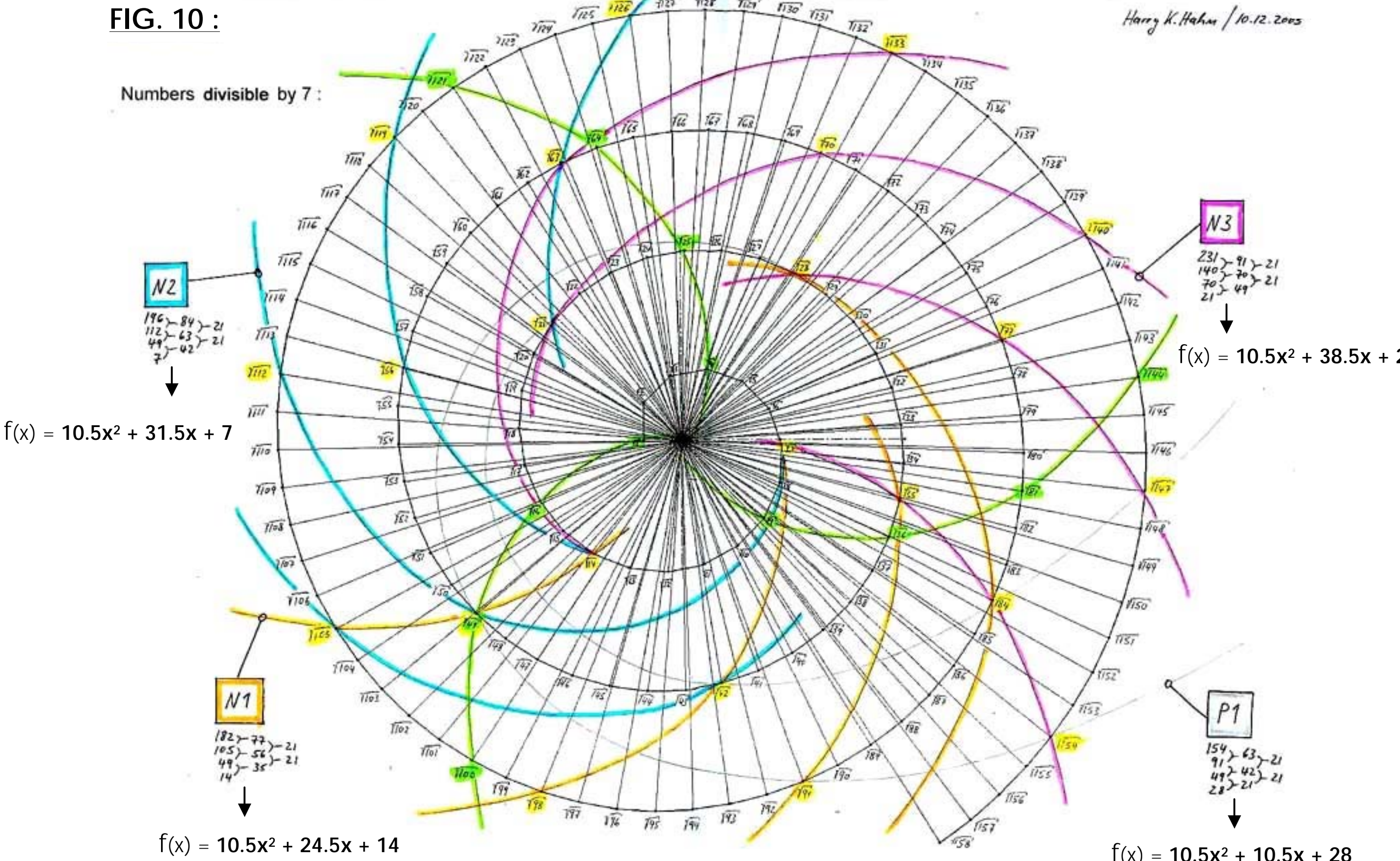

### 4.21 Properties of the spiral graph systems containing the numbers divisible by 7 :

- **The numbers divisible by 7** lie on **3** spiral-graph-systems with a **negative** rotation direction ( drawn in orange, blue and pink ) and on **3** spiral-graph-systems with a **positive** rotation direction ( only one system drawn in light grey lines ! ).
  The **3** negative spiral graph systems are named **N1, N2** and **N3** and the **3** positive systems are named **P1, P2** and **P3** ( only two spiral arms of P1 are drawn ! ).
  → Note : Not all spiralarms of the described spiral-graph-systems are drawn !
- The spiral-graph-systems **N1, N2** and **N3** as well as **P1, P2** and **P3** are arranged in an angle of around 120° to each other ( in reference to the centre of the square root spiral ) , and they seem to refer to the three-symmetrical arrangement of the **3** spiral-graphs of the square numbers 4, 9, 16, 25, 36, 49,… ( drawn in green ).
- Calculating the differences of the consecutive numbers lying on one of the drawn spiral arms, and then further calculating the differences of these differences ( **2. Differential** ), results in the constant value **21** for the positive as well as the negative rotating spiral-graphs.
  → see difference calculation for **4 exemplary spiral arms** ( one spiral arm per system ) in FIG. 10 beside the names of the spiral-graph-systems N1, N2, N3 and P1. ( → P2 & P3-system not shown ! )
  These **4** exemplary spiral arms are defined by the following quadratic polynomials :
  **Quadratic Polynomials of exemplary spiralarms :**

  $10.5x^2 + 24.5x + 14$ → belongs to **N1** - system
  $10.5x^2 + 31.5x + 7$ → belongs to **N2** - system
  $10.5x^2 + 38.5x + 21$ → belongs to **N3** - system
  $10.5x^2 + 10.5x + 28$ → belongs to **P1** - system

Natural numbers divisible by a certain prime factor are not distributed in a random way across the square root spiral ! This is evident from FIG. 10 !
The arrangement of the natural numbers divisible by 7 is a good example which shows the highly symmetrical distribution of certain number groups across the square root spiral in defined spiral systems.
In this case it is a highly three-symmetrical distribution similar to the distribution of the square numbers contained in the three spiral-graphs Q1, Q2 and Q3 drawn in green !

### 4.3 The distribution of natural numbers divisible by the prime factors 13 and 17 :

The next two images show the analysis results regarding the distribution of the natural numbers which are divisible by the prime factors 13 and 17 : ( → **only diagrams shown !** )

**FIG. 11 :**

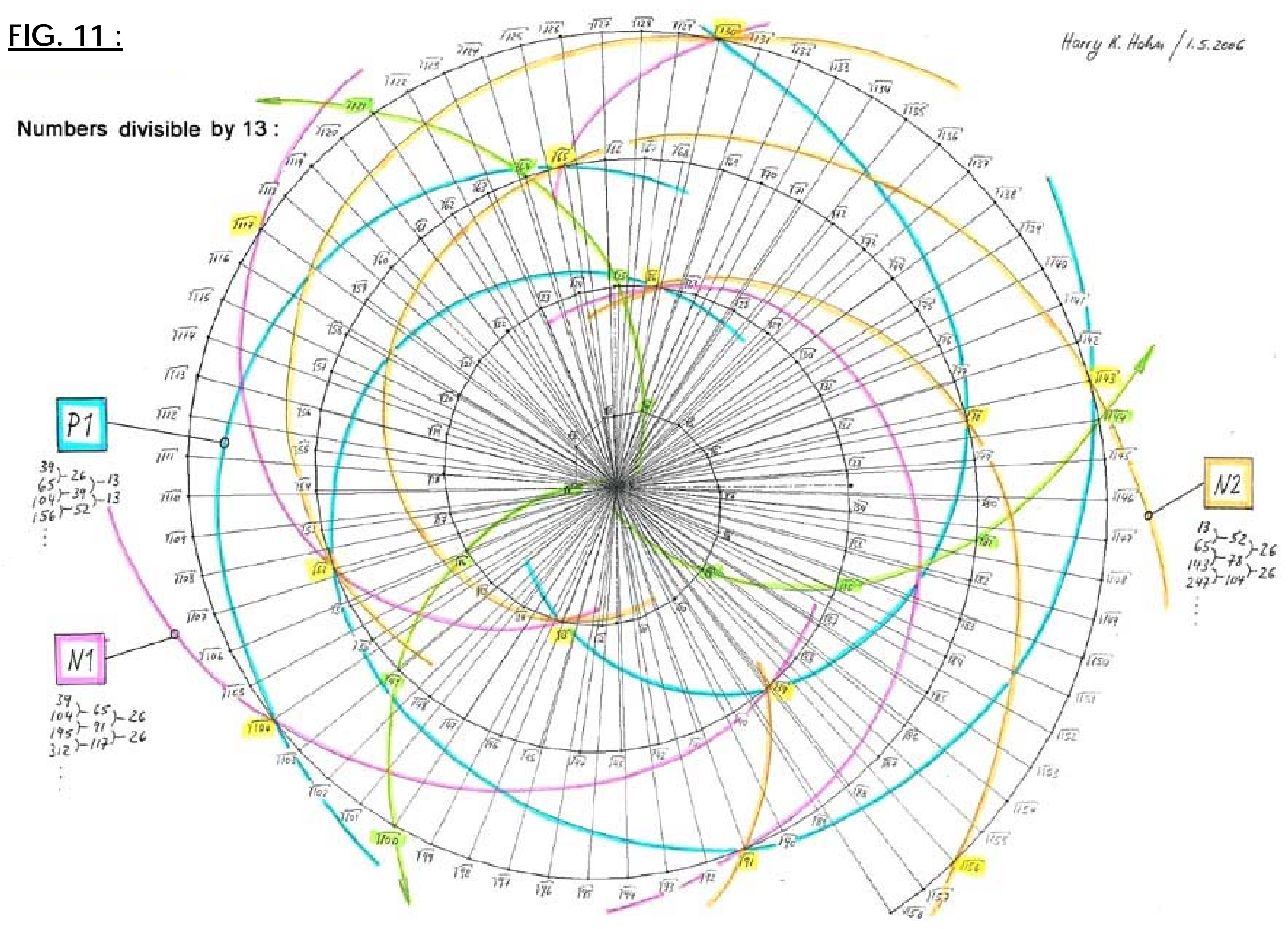

**FIG. 12 :**

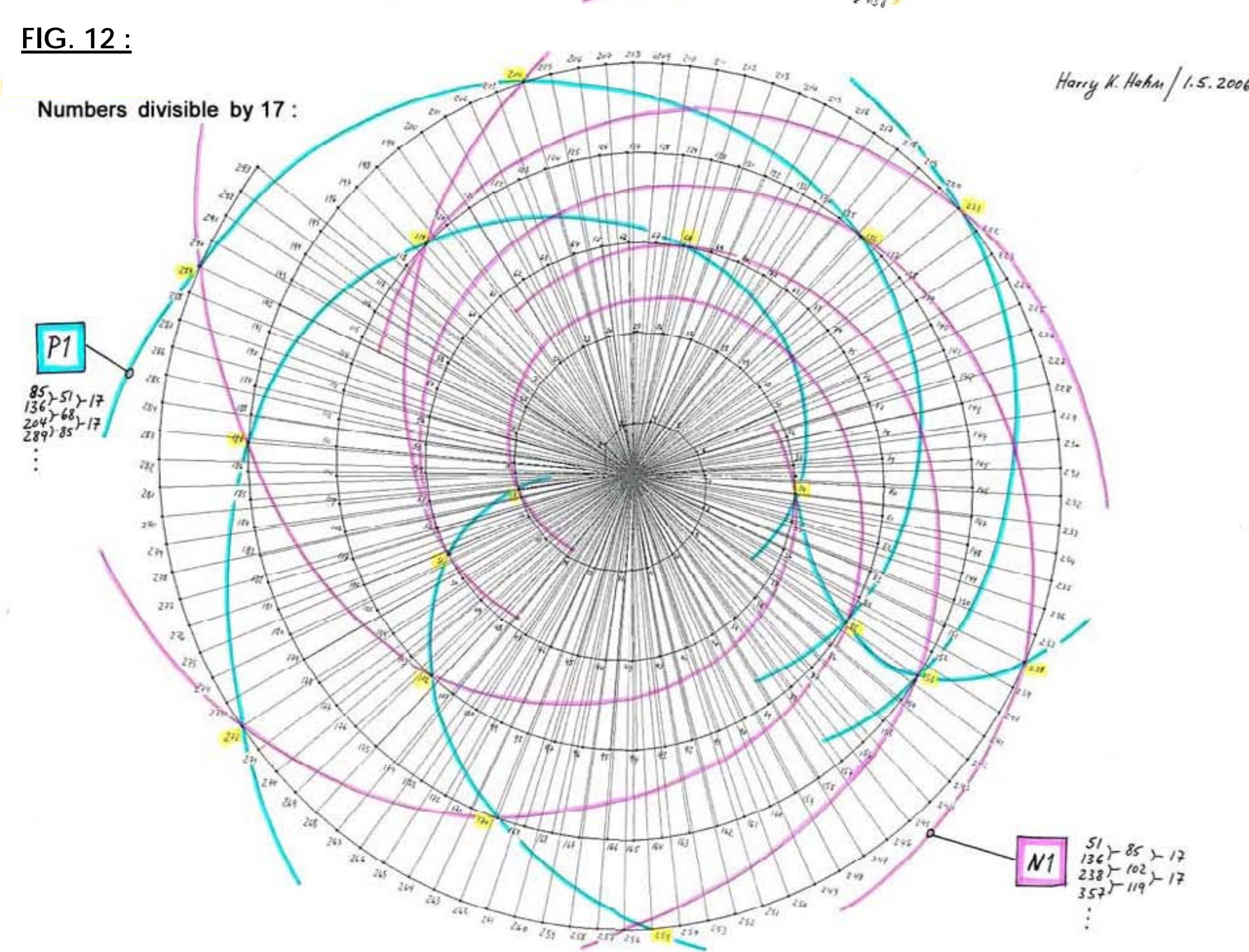

### 4.4 The distribution of natural numbers divisible by the prime factors 2, 3, 5, 13 and 17 :

In the same way as shown in FIG. 9 and FIG.10 for the natural numbers divisible by 11 or 7, I have also carried out a detailed analysis for the numbers divisible by the prime factors 2, 3, 5, 13 and 17.

These detailed analyses together with **high resolution images** of the spiral graph systems can be found in the **arXiv – archive** under my author name .

## 5 What causes the described Spiral Graph Systems ?

→ " The spiral graphs shown in FIG. 1 / 9 / 10 / 11 and 12 are caused by quadratic polynomials.
In principle every quadratic polynomial causes a sequence of radii, which takes an archimedian spiral-like course, when marked on the Square Root Spiral !

And the spiral angle of this so-created spiral graph converges ! "

This is essentially the conclusion of the mathematical analysis.

**→ explanation see mathematical section : Chapter 2 " The Spiral Arms"**

The exact course of these quadratic polynomials is given by the structure of the Square Root Spiral.

To better understand the whole structure of the Square Root Spiral, the following graph can be used :

The " **difference graph to the X-axis** " ( → see **FIG. 16** at page 30 in the **Appendix,** ), shows that the length of the circumference of the Square Root Spiral actually increases by approximately 20 per wind of the spiral. → see " **2. Differential** " **on FIG. 16**. The "difference graph to the X-axis" describes the difference of the increase of the circumference of the square root spiral to the number 20 per wind, in reference to the x-axis of the graph. This special graph, which represents the quadratic polynomial $f(x) = 10x^2 - 14x + 6$ can be used for a further analysis of the structure and the growth-behavior of the Square Root Spiral.

As already mentioned in the description of the highly three-symmetrical spiral graphs Q1, Q2 and Q3, which contain the square roots of the quadratic numbers, there must be a profound logic which governs the structure of the Square-Root-Spiral ! And this logic is definitely not understood yet !

And as a proof for this assumption, I can show a general rule which governs the existence of the described spiral graph systems as shown in FIG 9 to 12 !
These spiral-graph-systems shall be called "Number-Group-Spiral-Systems" to indicate that these spiral graph systems represent certain number groups ( e.g. numbers divisible by 2, 3, 5, 7, 11, 13, 17,… )

Before I explain this rule, I want to emphasize that this rule highly depends on the value of the **2. Differential** of the numbers lying on these spiral graphs as shown in the examples in FIG. 9 to 12

As described for the spiral graphs in FIG. 9 and FIG. 10, the " **2. Differential** " of all spiral graphs is constant and it is equal for all spiral graphs ( quadratic polynomials ) with the same rotation direction.
The " 2. Differential " can easily been calculated by calculating the differences of the consecutive numbers lying on one of the drawn spiral arms, and then further calculating the differences of these differences.

→ See example in FIG. 10 beside the name of the spiral-graph-sytem N1. In this example the differences of the numbers 182, 105, 49 and 14 on this spiral arm are 77, 56 and 35. And the differences of these numbers are all **21**. So the 2. Differential of this spiral-graph system is **21**.

With the help of the Newton Interpolation Polynomial and the calculated first and second differences the quadratic polynomials belonging to these spiral graphs can then be calculated.

### 5.1 Overall view of the distribution of the natural numbers on the Square Root Spiral

**Table 2** on page 32 in the **Appendix** shows the analysis results referring to the distribution of certain number groups on defined spiral-graph-systems on the "Square Root Spiral" ( for example the distribution of the square numbers or the distribution of natural numbers divisible by the prime factor 11 etc. )
It also shows the number sequences of exemplary spiral arms of the found spiral-graph-systems.
( → number sequences of one spiral arm per spiral-graph-system )

**Table 3A & 3B** at page 33/34 in the **Appendix** shows the quadratic polynomials of the exemplary spiral arms shown in **Table 2**.

This allows a first overall view of the quadratic polynomials, which define the spiral arms in the spiral-graph-systems shown in FIG. 9 to 12

## 5.2 The general rule which determines the existence of the described Spiral Graph Systems

→ There is an interdependency between the number of spiral-graph-systems with the same rotation direction , for a certain number group ( → e.g. all numbers divisible by the prime factor 11 ) and the " **2. Differential**" belonging to these spiral-graph-systems.

**<u>This interdependency can be expressed by the general formula shown in the head of the following table :</u>**

The table below clearly shows, that the discribed interdependency applies for all analysed number groups ( numbers divisible by the prime factors 2, 3, 5, 7, 11, 13 and 17 ) :

| **prime factor of number group** | **X** | **number of spiral graph systems** [ with a **negative (N)** or a **positive (P)** rotation direction ] | | **=** | **" 2. Differential "** |
|---|---|---|---|---|---|
| 2 | X | 10 | (N) | = | 20 |
| 2 | X | 9 | (P) | = | 18 |
| 3 | X | 7 | (N) | = | 21 |
| 3 | X | 6 | (P) | = | 18 |
| 5 | X | 4 | (N or P) | = | 20 |
| 7 | X | 3 | (N or P) | = | 21 |
| 11 | X | 2 | (N or P) | = | 22 |
| 13 | X | 2 | (N) | = | 26 |
| 13 | X | 1 | (P) | = | 13 |
| 17 | X | 1 | (N or P) | = | 17 |
| 19 | X | 1 | (N or P) | = | 19 |

( X = multiplication symbol )

Beside this general rule which determines the number of spiral-graph-systems , there is also a mathematical explanation available, which describes the character of single spiral-graphs.
**→ see Mathematical Section " The Spiral Arms ".**

Further there are also some notable differences in the number of spiral-graph-systems with an opposite rotation direction per number group :

## 5.3 Listing of differences of the number of spiral-graph-systems with an opposite rotation direction per number group :

- **The natural numbers divisible by 5, 7, 11, 17** and **19** lie on the same number of spiral-graph-systems for the **negative** as well as for the **positive** rotation direction. This also seems to be the case for all natural numbers divisible by prime factors >**19**.
- **The natural numbers divisible by 2, 3** and **13** lie on a different number of spiral-graph-systems for the **negative** and the **positive** rotation direction ( for example natural numbers divisible by the prime factor **3** lie on **7** spiral graph systems with a negative rotation direction and on **6** spiral graph systems with a positive rotation direction )
- **The natural numbers divisible by 2, 3, 5, 7** and **11** lie on more than one spiral graph system with either a positive rotation direction or a negative rotation direction.
- **The natural numbers divisible by prime factors ≥ 17** lie on only one spiral graph system for either the positive rotation direction or the negative rotation direction. (→ see example FIG.12 )

Also interesting is the fact that the " 2.Differential" of the spiral-graph-systems belonging to numbers divisible by a prime factor ≥ 17 is equal to the prime factor itself.
( e.g. natural numbers divisible by the prime factor 17 lie on one positive and one negative rotating spiral graph system with the constant number 17 as the " 2.Differential" of the graphs )

## 6 The distribution of Prime Numbers on the Square Root Spiral :

The distribution of the prime numbers on the square root spiral should interest every professional mathematician working in the field of number theory !!

Prime Numbers clearly accumulate on spiral graphs, which run through the square root spiral ( Einstein Spiral ) , in a similar fashion as the square numbers, or natural numbers which are divisible by the same prime factors ( as shown in FIG. 9 to 12 ).

**And all these "Prime Number Spiral Graphs" represent quadratic polynomials with special coefficients.**

Because I have described the distribution of prime numbers on the Square Root Spiral in more detail in another paper, I only want to show here one type of spiral-graph-system, which describes the distribution of the prime numbers on the Square Root Spiral.
It is probably the most impressive one. But there are other such systems existing with a different value of the " **2. Differential** " !

My **complete analysis of the Prime Number Spiral Graphs** can be found **on the arXiv–archive** under my author name and the following title :

**→ " The Ordered Distribution of Prime Numbers on the Square Root Spiral "**

The following picture **FIG.13** shows how the prime numbers are clearly distributed on defined spiral-graph systems which are arranged in a highly symmetrical manner around the centre of the Square Root Spiral.

On the shown **3 Prime-Number-Spiral-Systems** ( **PNS** ) **P18-A, P18-C** and **P18-C** , the prime numbers are located on pairs of spiral arms, which are separated by **3** numbers in between. And two spiral arms of one such pair of spiral arms, are separated by **1** number in between.

All spiral-graphs of the shown **3** Prime-Number-Spiral-Systems ( **PNS** ) have a **positive rotation direction** ( **P** ) and the " **2. Differential** " of all spiral-graphs is **18**.
That's why the first part of the naming of the **3** Prime-Number-Spiral-Systems ( **PNS** ) is **P18**.
The **3** spiral-graph-systems **A** ( drawn in orange ), **B** ( drawn in pink ) and **C** ( drawn in blue ) have further spiral arms. But for clearity there are only around 10 spiral arms drawn per system.

One striking property of all spiral arms is the nonexistence of numbers which are divisible by **2** or **3**.

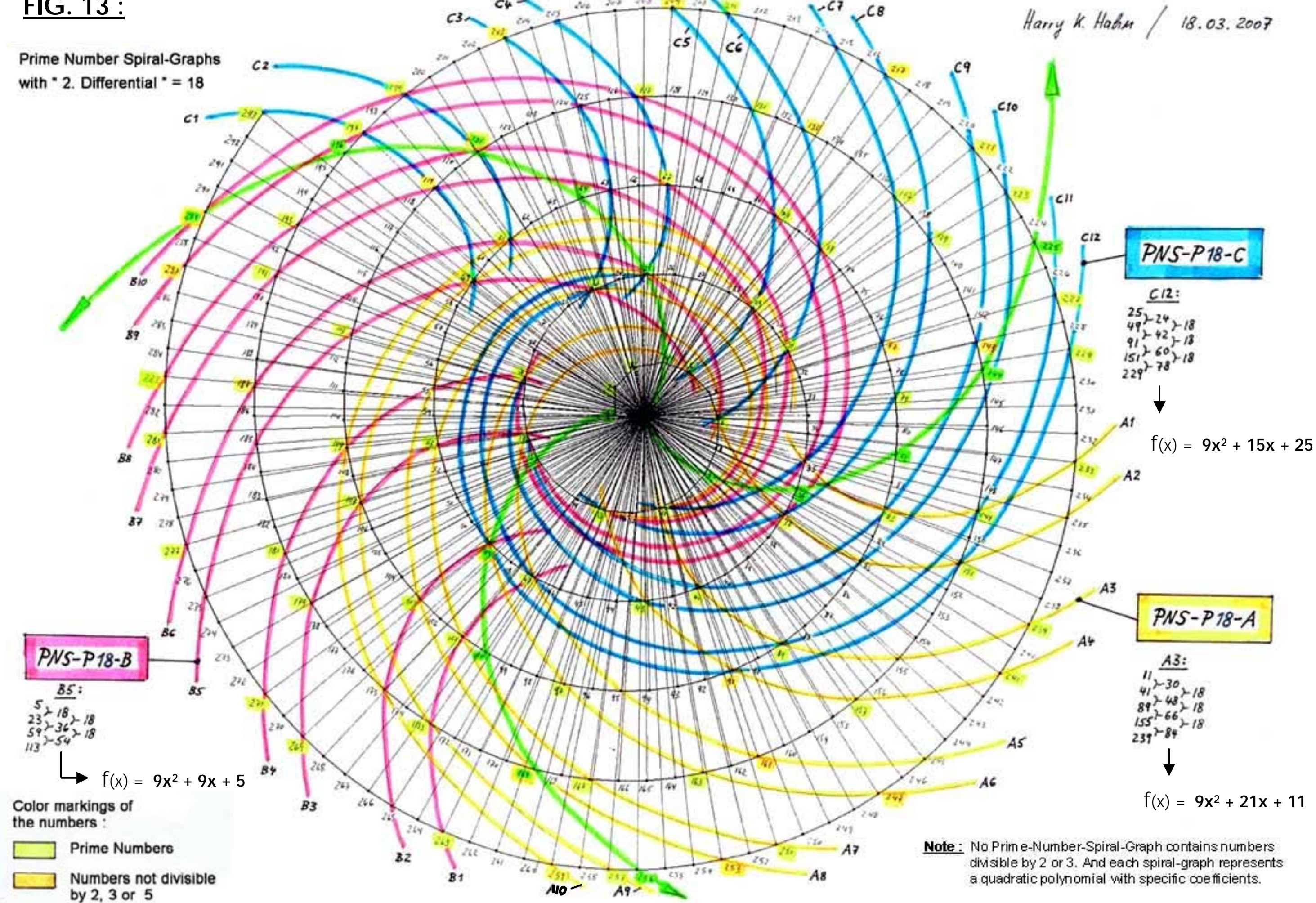

## 7 The distribution of Fibonacci Numbers on the Square Root Spiral :

Before I give you an insight into the relationship between the Square Root Spiral and Fibonacci Numbers, I want to invite you to read my detailed analysis about Fibonacci Numbers. This analysis contains some interesting facts and some new discoveries about Fibonacci Numbers !!

I intend to file this analysis in the arXiv–archive under the following title :

→ **" The mathematical origin of natural Fibonacci Sequences, and the periodic distribution of prime factors in these sequences."**

Fibonacci number sequences seem to play an important role in the structure of the Square Root Spiral.

Fibonacci numbers divide the Square Root Spiral into areas whose proportions strive for a constant ratio for $\sqrt{x} \rightarrow \infty$

And the ratio of the angles of two such successive areas on the Square Root Spiral is striving for a constant number at infinity too ! ( → see explanation below )

In both cases, this ratio is closely connected with the "golden mean" ( or golden section, or golden ratio )

The occurrence of these ratios, indicates that there is a special relationship between the Square Root Spiral and the Fibonacci Sequences !

**Ratio of the angles between consecutive Fibonacci Numbers → Self-Avoiding-Walk-Constant " SAW-F1" :**

If we mark the square roots of the Fibonacci Numbers 1, 2, 3, 5, 8, 13, 21,... on the Square Root Spiral and then measure the angles between the square roots of the numbers 1 and 2, 2 and 3, 3 and 5, 5 and 8, 8 and 13, 13 and 21 ...and so on, then we will get the following angles as a result :

$\alpha_1$ = 45° ; $\alpha_2$ = 35,26° ; $\alpha_3$ = 56,57° ; $\alpha_4$ = 67,01° ; $\alpha_5$ = 88,34° ; $\alpha_6$ = 111,40° etc.

**FIG. 14-A**

Calculated angle ratios :

$$\frac{\alpha_2}{\alpha_1} = 0{,}784$$

$$\frac{\alpha_3}{\alpha_2} = 1{,}604$$

$$\frac{\alpha_4}{\alpha_3} = 1{,}185$$

$$\frac{\alpha_5}{\alpha_4} = 1{,}318$$

$$\frac{\alpha_6}{\alpha_5} = 1{,}261$$

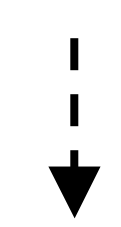

If we now calculate the ratios of successive angles, as shown above, we get the following ratios as result :

0,784 ; 1,604 ; 1,185 : 1,318 ; 1,261 etc

It is easy to see that this ratio is quickly approaching a constant number for $\sqrt{x} \rightarrow \infty$

The square root spiral which I used for my analysis is precisely constructed up to sqrt293 and allows one to calculate this constant with the following accuracy :

**1,272242 < SAW-F1 < 1,272507**

I call this constant "SAW-F1" which means "self-avoiding-walk constant " F1
( Here F1 stands for Fibonacci Sequence 1 ).

This constant is already known as " self-avoiding-walk-constant 1.272…"

→ see book : Mathematical Constants" from Steven R. Finch

There might still be a bit of inaccuracy left in the calculated range of **1,272242 < SAW-F1 < 1,272507** for the angle ratio, because of inaccuracies in the sums of the angles, calculated by the CAD-System.

So the real value of this constant **SAW-F1** could finally be slightly higher or lower.

It could turn out, that the true value of this constant **SAW-F1** is **1.27201965…** !
This number is the square root of the golden mean ( golden section ) !

→ The golden mean ( golden section ) : $\tau$ **= 1,61803399...** and $\sqrt{\tau}$ **= 1.27201965.....**

There is a good reason that this could finally be the correct constant ! → see next paragraph

**Ratio of the areas between consecutive Fibonacci Numbers → Area-Ratio-Constant "ARC-F1" :**

My reason for the above-mentioned assumption is the value of the constant for the proportions of the areas between the square roots of the Fibonacci Numbers. Because the proportions of the areas also strive for a constant which also seems to be linked to the "golden ratio" !

We mark again the square roots of the Fibonacci Numbers 1, 2, 3, 5, 8, 13, 21,… on the square root spiral.
Then we calculate the areas between these marked square roots ( areas marked in red, green, blue etc.)

And if we now calculate the ratios between successive marked areas as shown below, then we get the following ratios as a result : 1.4114 , 2.6389 , 1.9644 , 2.1512 , 2.055 ….

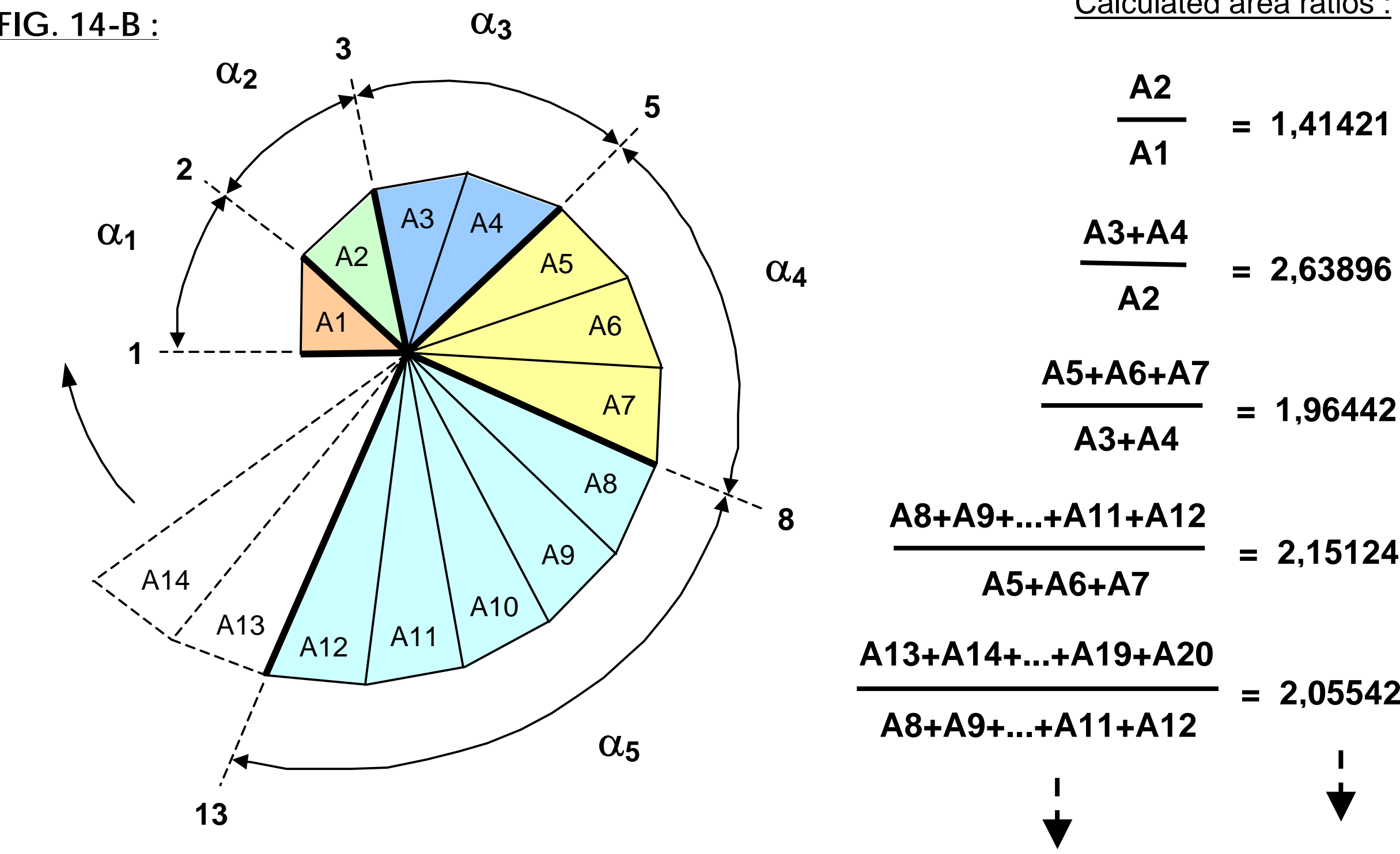

The ratio between successive areas ( as shown above ) is striving for a constant number for $\sqrt{x} \rightarrow \infty$

I calculated this constant with the following accuracy : **ARC-F1 = 2.05819 ± 0.0003**

I call this constant "ARC-F1" which means "area ratio constant" -F1
( Here F1 again stands for Fibonacci Sequence 1 ).

As mentioned before : This constant is closely related to the golden mean ( golden section ) !

With a high probability it is equal to $\tau \times \sqrt{\tau}$ = 2.058171… with $\tau$ **= 1,61803399..** ( = golden mean )

There is an interesting study existing to this constant $\tau \times \sqrt{\tau}$ from Mark A. Reynolds, which has the following title : **"The unknown Modulor: the "2.058" Rectangle**

This study can be downloaded from this weblink : **http://www.springerlink.com/content/w534664pmjx0/**

# 8 Final comment / References

The square root spiral ( or Einstein Spiral ) shows the interdependencies between natural numbers in a visual way. Therefore, it can be considered to be a kind of visual representation of number theory !

Through pure graphical analysis of this amazing structure, the higher logic of the ( spatial ) distribution of natural numbers ( and special sub-groups like square numbers or prime numbers ) comes to light and is very easy to grasp, because it is visible !!

That's why mathematicians who read this paper should continue my work and do a more extensive analysis of the square root spiral, by using a precise computer model of the square root spiral and by using more advanced analysis software and analysis techniques, in a similar way as I have !

I haven't found any scientific study which carried out a similar graphical analysis of the ( spatial ) distribution of natural numbers on the square root spiral as shown here in my study. This offers a great opportunity because there is a lot of unknown land to discover here ! And I probably have only just reached the beach of this new land and made a first clumsy step !

In December 2005 and June 2006 I sent this study ( in CD-format ) to a dozen universities in Germany for an assessment. But there wasn't much response ! That's why I decided to publish my discoveries here.

Prof. S.J. Patterson from the University of Goettingen found my discoveries very interesting. Prof. S.J. Patterson was especially interested in the spiral graphs which contain the Prime Numbers. These spiral graphs are special quadratic polynomials, which are of great interest to Prime Number Theory.
For example the quadratic polynomial B3 in FIG. 15-D → B3 = $F(x) = 9x^2 + 27x + 17$ ( or $9x^2 + 9x - 1$ )
or the quadratic polynomial K5 in FIG. 15-F → K5 = $F(x) = 11x^2 + 25x + 13$ ( or $11x^2 + 3x - 1$ )
→ see my study " The Ordered Distribution of Prime Numbers on the Square Root Spiral "

Prof. Ernst Wilhelm Zink from the Humboldt-University in Berlin also found my study very interesting and he organized a mathematical analysis of the spiral-graphs shown in this study as well as an analysis of some propositions which are described in my study referring to the Fibonacci Number Sequences. This mathematical analysis was carried out by Mr. Kay Schoenberger a student of mathematics on the Humboldt-University of Berlin. Mr. Kay Schoenberger is currently doing his mathematical dissertation.

On this occasion I want to thank Mr. Kay Schoenberger for the smart proof of my propositions and Prof. Ernst Wilhelm Zink for his interest in my work, for all his help to organize the mathematical analysis, which considerably upgraded the value of this study, and for his support during the publication process.

## References

In the following I have listed some literature and weblinks, which refer to other studies of the Square Root Spiral, except of the first reference, which is not related to the Square Root Spiral.
The first reference is a very interesting website which deals with a special "Number Spiral", which is related to the well known Ulam Spiral. Everyone who reads my analysis should definitely have a closer look to this website too, because many of the analysis results shown on this website are in connection with my findings !

( In comparison to the Square Root Spiral, the mentioned "Number Spiral" is just wound 3 times tighter ! )

# Mathematical Section

## 1 The Correlation to $\pi$

In the following part I would like to discuss in detail two assertions made by Mr. Hahn in the brief description. He claims here that

- the distance between two successive winds of the Square Root spiral converges to $\pi$,
- the angle between two successive integer cathets converges to 2.

Using relatively elementary methods, both assertions are proven below.

If $t_n$ is the angle of the *nth* spiral segment (in radian), we get

$$\tan(t_n) = \frac{1}{\sqrt{n}}$$

(ratio of counter cathet/ancathet ). From this it follows that

$$t_n = \arctan\left(\frac{1}{\sqrt{n}}\right)$$

and the angle *w(k)* between the first and the *(k+1)th* cathet is

$$w(k) = \sum_{n=1}^{k} t_n$$

As $t_n$ represents a non increasing sequence which tends to zero, the integral approximation

$$w(k) = \sum_{n=1}^{k} t_n = \int_0^k \arctan\left(\frac{1}{\sqrt{n}}\right) \mathrm{d}n + c_1(k)$$

applies, whereby for the remainder term $\lim_{k\to\infty} c_1(k) = \text{const.}$ is valid.

The integral can be solved and from this we get

$$w(k) = -\arctan(\sqrt{k}) + k \cdot \arctan\left(\frac{1}{\sqrt{k}}\right) + \sqrt{k} + c_1(k)$$

As for $x > 0$ the equality $\arctan(x) = \frac{\pi}{2} - \arctan(\frac{1}{x})$ applies, we can also write

$$w(k) = \arctan(\tfrac{1}{\sqrt{k}}) - \tfrac{\pi}{2} + k \cdot \arctan\left(\frac{1}{\sqrt{k}}\right) + \sqrt{k} + c_1(k)$$

The next step is to represent *k* as a function of *w*. This certainly does not work directly but we can rely on the Taylor series

$$\arctan(x) = x - \frac{x^3}{3} + \frac{x^5}{5} - + \ldots$$

Then we get the following equality

$$\arctan\left(\frac{1}{\sqrt{k}}\right) = \frac{1}{\sqrt{k}} + \mathcal{O}\left(\frac{1}{k\sqrt{k}}\right)$$

Substitute:

$$\begin{aligned} w(k) &= \underbrace{\frac{1}{\sqrt{k}} + \mathcal{O}\left(\frac{1}{k\sqrt{k}}\right)}_{\to 0} - \frac{\pi}{2} + \frac{k}{\sqrt{k}} + \underbrace{\mathcal{O}\left(\frac{1}{\sqrt{k}}\right)}_{\to 0} + \sqrt{k} + c_1(k) \\ &= 2\sqrt{k} + c_2(k) \quad (\text{ where } \lim_{k\to\infty} c_2(k) = \text{const.}) \\ k(w) &= \frac{1}{4}(w - c_3(w))^2 \quad (\text{ where } \lim_{w\to\infty} c_3(w) = \text{const.}) \end{aligned}$$

For the radius $r$ ( i.e. the big cathet) we have $r(k) = \sqrt{k}$, thus

$$r(k(w)) = r(w) = \sqrt{\frac{1}{4}(w - c_3(w))^2} = \frac{1}{2}(w - c_3(w))$$

The radius therefore is proportional to the angle (the *Archimedean Spiral* also has this property).
The distance *a(w)* of the spiral arms is the difference of the radii after a full rotation, therefore

$$\begin{aligned} a(w) &= r(w + 2\pi) - r(w) \\ &= \frac{1}{2}\left(w + 2\pi - c_3(w + 2\pi)\right) - \frac{1}{2}(w - c_3(w)) \\ &= \pi - \underbrace{\frac{1}{2}\left(c_3(w + 2\pi) - c_3(w)\right)}_{\to 0} \end{aligned}$$

Hence we obtain: $\lim\limits_{w\to\infty} a(w) = \pi$

Accordingly for the angle $v(k)$ between the square numbers $k^2$ and $(k+1)^2$ the following applies:

$$\begin{aligned} v(k) &= w((k+1)^2) - w(k^2) \\ &= 2(k+1) + c_2((k+1)^2) - 2k - c_2(k^2) \\ &= 2 + \underbrace{c_2((k+1)^2) - c_2(k^2)}_{\to 0} \end{aligned}$$

The limit is therefore $$\lim_{k \to \infty} v(k) = 2 = \tfrac{2\pi}{\pi} \equiv \tfrac{360^\circ}{\pi}$$

This therefore corresponds exactly with Mr. Hahn's prediction.

I presume that by the not completely clearly defined angle between the integer cathets of two successive spiral arms Mr. Hahn probably means the angle between two radii whose length differs by exactly *3*. On account of $\pi \approx 3$, these form a very small angle ( not paying attention to the full rotation in between).

From the above analysis it immediately follows that this converges to $2\pi - 6$ or

$$360^\circ - 3 \cdot \tfrac{360^\circ}{\pi}$$

as described by Mr. Hahn.

## 2 The Spiral Arms

The pictures FIG. 1 and FIG. 9 to 12 show a graphical analysis of the Square Root Spiral carried out by Mr. Hahn.
In FIG. 1 Mr. Hahn marked all Square Numbers ( green ) and divided them into three groups, through which he draw three spiral shaped graphs ( drawn in green lines ).
He made a similar subdivision for other number groups in FIG. 9 to 12.
In FIG.9 for example he divided all numbers, which are divisible by 11, into two systems of spiral graphs ( drawn in orange and pink ).

From the previous part " The Correlation to $\pi$", it follows that two Square Numbers $k^2$ and $(k + 3)^2$ enclose a small angle $\approx 2\pi - 6$

Mr. Hahn therefore simply subdivides the set of integer radii into three equivalence classes, using the classes of $\mathbb{Z}/3\mathbb{Z}$ .

The set of integer radii can also be partitioned using $\mathbb{Z}/n\mathbb{Z}$, where $n$ is a natural number. But this doesn't result in similar " nice " angles.

Because the angle behaves nearly propotional to the radius, we therefore always get similar spiral arms.

The spiral graphs which contain the numbers divisible by 11 can be analysed in the same way. Mr. Hahn obviously used the following procedure to construct these spiral graphs :

He started with a multiple of 11 and then located the closest successor number on the next wind of the Square Root Spiral. Coming from the center of the Square Root Spiral the closest successor number which is divisible by 11 would lie a bit on the right ( for the N-class spiral arms ), or a bit on the left ( for the P-class spiral arms ), on the next wind of the Square Root Spiral, in reference to the start number. Mr. Hahn connected these numbers with a graph.
By continuing this graph the further successor numbers can then be located.
The continuing of this graph is made in such a way that nearly a constant spiral angle is achieved. That there are always corresponding numbers divisible by 11 located on the Square Root Spiral, through which the graph can be continued, is explained below.

In this way a sequence of multiples of 11 is created, which all lie on a defined spiral arm. For example graph $N_1$ drawn in FIG. 9 results in the following number sequence :

$$N_1 : 22, 77, 154, 253, ...$$

Mr. Hahn now sets up the following sequences of differences

| Sequence $N_1$ : | 22 | | 77 | | 154 | | 253 | | ... |
|---|---|---|---|---|---|---|---|---|---|
| first difference : | | 55 | | 77 | | 99 | | ... | |
| second difference | | | 22 | | 22 | | ... | | |

,

and he noticed that all second differences are equal to 22. And it is the same with the other sequences drawn in FIG. 9.
To explain this fact we go back to the following expression :

$$k(w) = \tfrac{1}{4}(w - c_3(w))^2$$

( → see part 1 " The Correlation to $\pi$ " )

The numbers $k$ are prportional to the square of the angle $w$. This means that $N_1$ is probably the sequence of values of a quadratic polynomial at the natural numbers. The same holds for the other sequences drawn in FIG. 9 to 12.
A quadratic polynomial always defines a sequence, where the sequence of the second differences is constant.

For the proof we need some preliminary explanation: given arguments $t_1, t_2, ..., t_n$ and values $f(t_1), ... , f(t_n)$ the *Newton's divided differences* are defined as follows:

$$f[t_{i_1}, t_{i_2}] := \frac{f(t_{i_1}) - f(t_{i_2})}{t_{i_1} - t_{i_2}}$$

$$f[t_{i_1}, ..., t_{i_k}] := \frac{f[t_{i_1}, ..., t_{i_{k-1}}] - f[t_{i_2}, ..., t_{i_k}]}{t_{i_1} - t_{i_k}}$$

With their help the *Newton interpolation polynomial*

$$N(t) := f(t_1) + (t - t_1) f[t_1, t_2] + (t - t_1)(t - t_2) f[t_1, t_2, t_3] + ... \\ ... + \prod_{i=1}^{n-1} (t - t_i) f[t_1, ..., t_n]$$

can be defined. This fulfils the interpolation conditions

$$N(t_i) = f(t_i), \quad i = 1, ..., n$$

For the sequences given by Mr. Hahn, this yields

| | | | | | | | | | |
|---|---|---|---|---|---|---|---|---|---|
| sequence: | $f(1)$ | | $f(2)$ | | $f(3)$ | | $f(4)$ | | ... |
| first difference: | | $f[1,2]$ | | $f[2,3]$ | | $f[3,4]$ | | ... | |
| second difference: | | | $2f[1,2,3]$ | | $2f[2,3,4]$ | | ... | | |

The third differences are zero, so that the interpolation polynomial has no terms of third order (and this yields a quadratic polynomial). In general, equality of the *n-th* differences leads to a polynomial of degree $n$.

In our case we have with the short notation $f_i := f(t_i)$ the polynomial

$$\begin{aligned} N(t) &= f_1 + (t-1)f[1,2] + (t-1)(t-2)f[1,2,3] \\ &= f_1 + (t-1)(f_2 - f_1) + \tfrac{1}{2}(t-1)(t-2)(f_1 - 2f_2 + f_3) \end{aligned}$$

The generator polynomial for $N_1$ is therefore

$$\begin{aligned} N_1(t) &= 22 + (t-1)(77-22) + \tfrac{1}{2}(t-1)(t-2)(22 - 2\cdot 77 + 154) \\ &= 11(t^2 + 2t - 1) \end{aligned}$$

The similarity between difference and derivatives mentioned by Mr. Hahn is certainly visible.

If $t_1, \ldots, t_n \in [a,b]$ are pairwise different, and if $f \in C^{n-1}([a,b])$ , then according to a mean value theorem there exists a $\xi \in [a,b]$ ,

such that $$f[t_1, \ldots, t_n] = \frac{1}{(n-1)!} \cdot f^{(n-1)}(\xi)$$

## A motivation that Mr. Hahn's considered curves $N_i$ and $P_i$ are close to archimedian spirals?

For this we consider a quadratic polynomial.

$$p(t) = at^2 + bt + c$$

For the corresponding angle $w(t) := w(p(t))$ on the root spiral we have

$$w(t) = 2\sqrt{at^2 + bt + c} + c_2(p(t))$$

Now we consider the angle $v(t) := w(t+1) - w(t)$ between the radii of two successive sequence terms. For the limit we get the following.

$$\begin{aligned}
\lim_{t\to\infty} v(t) &= \lim_{t\to\infty} \left[ 2\sqrt{p(t+1)} - 2\sqrt{p(t)} + \underbrace{c_2(p(t+1)) - c_2(p(t))}_{\to 0} \right] \\
&= \lim_{t\to\infty} 2 \left[ \sqrt{p(t+1)} - \sqrt{p(t)} \right] \\
&= \lim_{t\to\infty} 2 \left[ \frac{p(t+1)-p(t)}{\sqrt{p(t+1)}+\sqrt{p(t)}} \right] \\
&= \lim_{t\to\infty} 2 \left[ \frac{2at+a+b}{\sqrt{at^2+2at+bt+b+c+1}+\sqrt{at^2+bt+c}} \right] \\
&= \lim_{t\to\infty} 2 \left[ \frac{t(2a+\mathcal{O}(\frac{1}{t}))}{\sqrt{t^2(a+\mathcal{O}(\frac{1}{t}))}+\sqrt{t^2(a+\mathcal{O}(\frac{1}{t}))}} \right] \\
&= \lim_{t\to\infty} 2 \left[ \frac{2a+\mathcal{O}(\frac{1}{t})}{2\sqrt{a+\mathcal{O}(\frac{1}{t})}} \right] \\
&= 2 \cdot \frac{2a}{2\sqrt{a}} = 2\sqrt{a}
\end{aligned}$$

The corresponding spiral angle therefore converges to $(2\sqrt{a} \mod 2\pi)$

**The number sequences given by Mr. Hahn are generated from quadratic polynomials !**

The mean value theorem qouted above implies immediately that the second differences are $2a$. For the example $N_1$ and that of the other sequences of multiples of $11$ , the parameter $a$ must according to the construction be a multiple of $\frac{11}{2}$ .

A corresponding investigation for the prime numbers 7, 13 and 17 is shown in Figures 10 to 12 as well as in tabular form in Table 2, 3A and 3B. Here the prime numbers give rise to a special choice of coefficients of the polynomial corresponding to them.

**Final Conclusion of the analysis of the spiral arms shown in FIG.1 and FIG. 9 to 12 :**

Every quadratic polynomial causes a sequence of radii, which takes an archimedian spiral-like course, when marked on the Square Root Spiral !
And the spiral angle of this so created spiral graph converges !

## 3 Area Equality

The following proposition refers to figure 7 in chapter 3.1.

**Proposition 1**

$$\lim_{M\to\infty} \frac{S(M+1)}{S(M)} = 1 \text{ for } S(M) := \sum_{n=M^2}^{M^2+2M} A_n$$

**Proof :** We consider the function

$$g(n) = A_n = \frac{1}{2}\sqrt{n}.$$

According to "Euler's summation formula" we have :

$$\sum_{n=1}^{N} g(n) = g(1) + \int_1^N g(x)dx + \int_1^N (x-[x])g'(x)dx.$$

Hence we get

$$(1) \qquad \sum_{n=M}^{N} \frac{1}{2}\sqrt{n} = \left[\frac{1}{3}x^{\frac{3}{2}}\right]_{M-1}^{N} + \frac{1}{4}\int_{M-1}^{N} (x-[x])x^{-\frac{1}{2}}dx,$$

and the error term is bounded from above by

$$\frac{1}{4}\int_{M-1}^{N} x^{-\frac{1}{2}}dx = \frac{1}{2}\left[x^{\frac{1}{2}}\right]_{M-1}^{N}.$$

Using the definition

$$T(M) := \sum_{n=M^2+1}^{(M+1)^2} A_n \;=\frac{1}{2}\sum_{n=M^2+1}^{(M+1)^2}\sqrt{n}$$

by $S(M) = T(M) - \frac{1}{2}$ we get

$$\lim_{M\to\infty} \frac{S(M+1)}{S(M)} = \lim_{M\to\infty} \frac{T(M+1)}{T(M)}.$$

Now we take (1) into account.

$$T(M) = \left[\frac{1}{3}x^{\frac{3}{2}}\right]_{M^2}^{(M+1)^2} + c = \frac{1}{3}((M+1)^3 - M^3) + c$$

The error is

$$c \leq \frac{1}{2}\left[x^{\frac{1}{2}}\right]_{M^2}^{(M+1)^2} = \frac{1}{2},$$

and the following equation finishes the proof.

$$\lim_{M\to\infty}\frac{T(M+1)}{T(M)} = \lim_{M\to\infty}\frac{(M+2)^3-(M+1)^3}{(M+1)^3-M^3} = \lim_{M\to\infty}\frac{3M^2+9M+7}{3M^2+3M+1} = 1$$

q.e.d.

# Appendix :

## FIG. 15 :

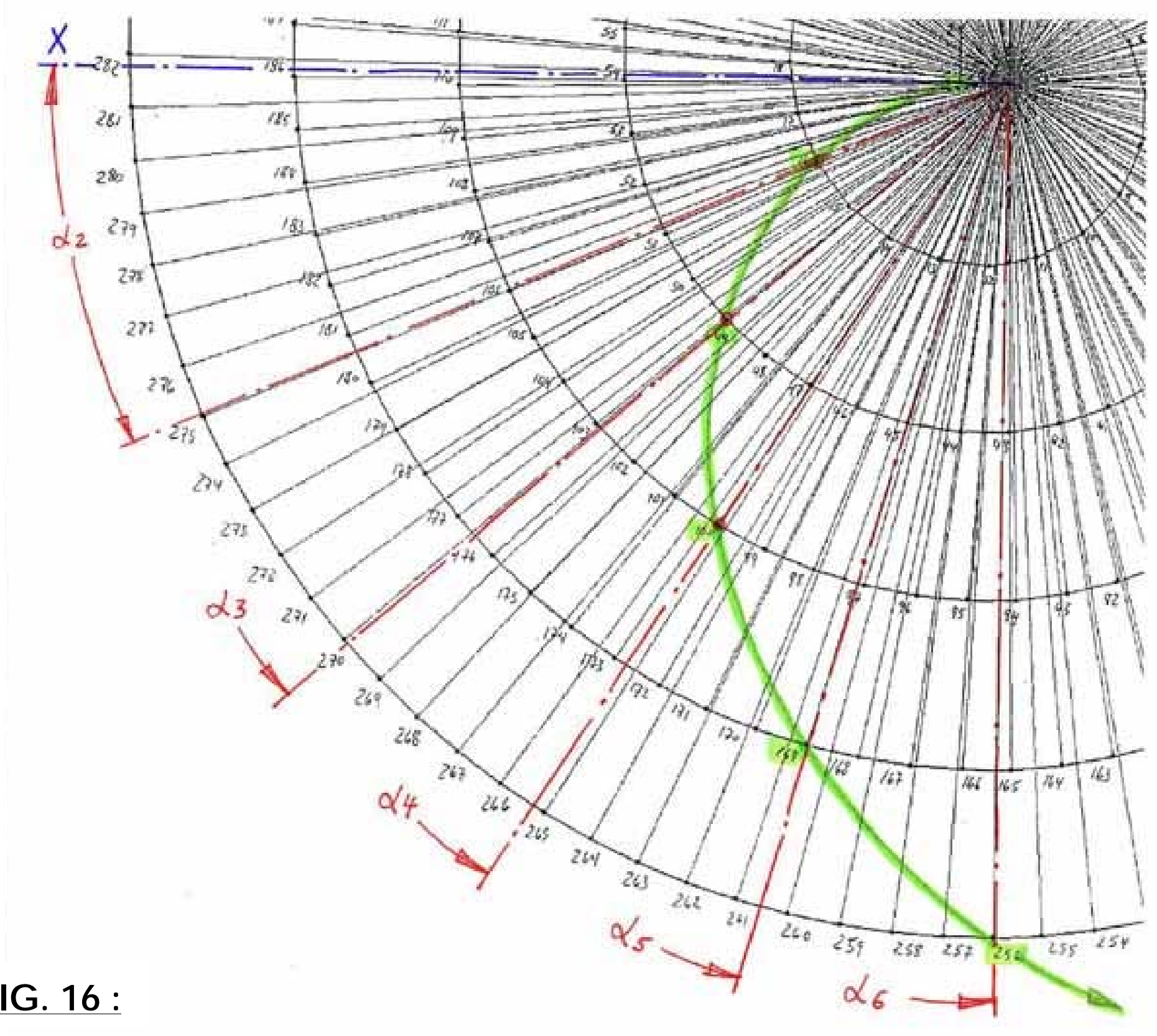

**Square-Numbers Spiral-Graph Q1:**

| Root | α | Angle to the X-Axis |
|---|---|---|
| $\sqrt{1}$ | 1 | 0° |
| $\sqrt{16}$ | 2 | 22,88582334° |
| $\sqrt{49}$ | 3 | 40,12815644° |
| $\sqrt{100}$ | 4 | 56,75743814° |
| $\sqrt{169}$ | 5 | 73,20458755° |
| $\sqrt{256}$ | 6 | 89,56818432° |

Calculated
quadratic polynomial
of Square-Number-Graph Q1 :

→ $f(x) = 9x^2 - 12x + 4$

## FIG. 16 :

**Difference-Graph to the X-Axis :**

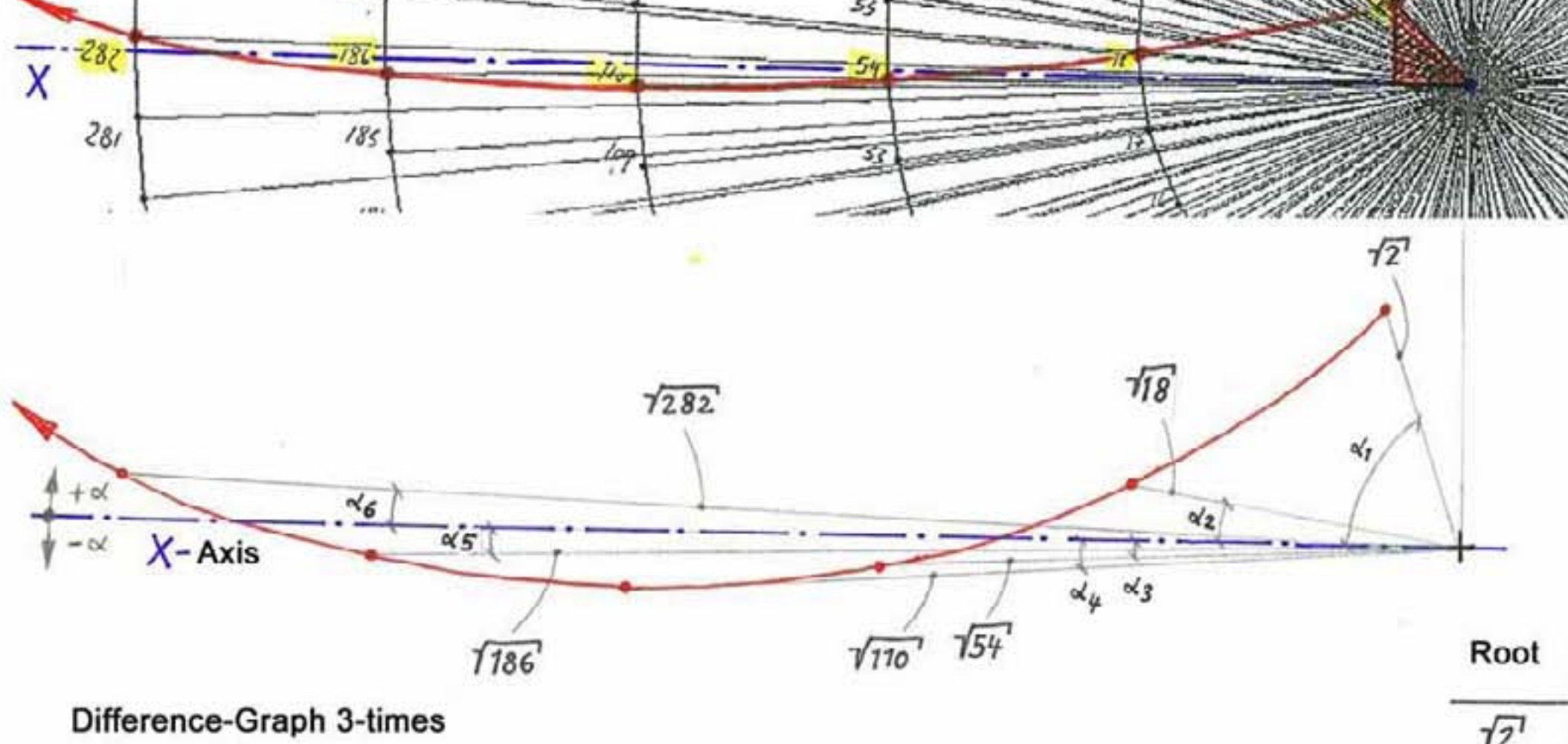

Difference-Graph 3-times
enhanced drawn :

The "2. Differential" of the numbers
on the Difference-Graph shows :

```
>>>   2
         }-  16
     18           }-  20
         }-  36
     54           }-  20
         }-  56
    110           }-  20
         }-  76
    186           }-  20
         }-  96
    282
```

| Root | α | Angle to the X-Axis |
|---|---|---|
| $\sqrt{2}$ | 1 | 45° |
| $\sqrt{18}$ | 2 | 4,78344235° |
| $\sqrt{54}$ | 3 | -0,26102520° |
| $\sqrt{110}$ | 4 | -0,87145670° |
| $\sqrt{186}$ | 5 | -0,10910812° |
| $\sqrt{282}$ | 6 | 1,25786779° |

Calculated quadratic polynomial
of Difference-Graph :

$f(x) = 2(5x^2 - 7x + 3) = \mathbf{10x^2 - 14x + 6}$

**Tabelle 1 :** Analysis to the development of the winding-distance of the „Einstein-Spiral“ ( Square Root Spiral )

| Length difference of two „root rays“ which differ by approximately one winding of the square root spiral | calculated winding-distance | winding No. : | calculated average winding-distance of this winding | accuracy of average winding-distance in comparison to $\pi$ |
|---|---|---|---|---|
| $\sqrt{21} - \sqrt{2}$ | 3,16836 | 2 | 3,1592037 | 99,44255 % |
| $\sqrt{24} - \sqrt{3}$ | 3,16693 | | | |
| $\sqrt{29} - \sqrt{5}$ | 3,14910 | | | |
| $\sqrt{38} - \sqrt{9}$ | 3,16441 | | | |
| $\sqrt{40} - \sqrt{10}$ | 3,16228 | | | |
| $\sqrt{42} - \sqrt{11}$ | 3,16412 | | | |
| $\sqrt{51} - \sqrt{16}$ | 3,14143 | | | |
| $\sqrt{53} - \sqrt{17}$ | 3,15700 | | | |
| $\sqrt{58} - \sqrt{20}$ | 3,14364 | 3 | 3,1443455 | 99,91245 % |
| $\sqrt{63} - \sqrt{23}$ | 3,14142 | | | |
| $\sqrt{68} - \sqrt{26}$ | 3,14719 | | | |
| $\sqrt{76} - \sqrt{31}$ | 3,15003 | | | |
| $\sqrt{79} - \sqrt{33}$ | 3,14363 | | | |
| $\sqrt{82} - \sqrt{35}$ | 3,13931 | | | |
| $\sqrt{97} - \sqrt{45}$ | 3,14065 | | | |
| $\sqrt{100} - \sqrt{47}$ | 3,14435 | | | |
| $\sqrt{103} - \sqrt{49}$ | 3,14889 | | | |
| $\sqrt{113} - \sqrt{56}$ | 3,14683 | 4 | 3,14428 | 99,91453 % |
| $\sqrt{127} - \sqrt{66}$ | 3,14539 | | | |
| $\sqrt{138} - \sqrt{74}$ | 3,14501 | | | |
| $\sqrt{142} - \sqrt{77}$ | 3,14141 | | | |
| $\sqrt{174} - \sqrt{101}$ | 3,14103 | | | |
| $\sqrt{178} - \sqrt{104}$ | 3,14363 | | | |
| $\sqrt{182} - \sqrt{107}$ | 3,14666 | | | |
| $\sqrt{191} - \sqrt{114}$ | 3,14320 | 5 | 3,142395 | 99,97447 % |
| $\sqrt{200} - \sqrt{121}$ | 3,14213 | | | |
| $\sqrt{209} - \sqrt{128}$ | 3,14312 | | | |
| $\sqrt{223} - \sqrt{139}$ | 3,14336 | | | |
| $\sqrt{228} - \sqrt{143}$ | 3,14141 | | | |
| $\sqrt{263} - \sqrt{171}$ | 3,14058 | | | |
| $\sqrt{268} - \sqrt{175}$ | 3,14195 | | | |
| $\sqrt{273} - \sqrt{179}$ | 3,14362 | | | |
| $\sqrt{284} - \sqrt{188}$ | 3,14099 | | | |
| $\sqrt{289} - \sqrt{192}$ | 3,14359 | | | |

**Note :** On the left side of this table, the length-difference of two „square root rays“ at a time are shown, which differ by approximately one winding of the „Square Root Spiral“ to each other ( → see FIG. 1 – Square Root Spiral ). With every further winding of the spiral, these length-differences ( which represent the winding-distance ), strive more and more for the constant $\pi$ . This is evident if we compare the average winding-distance with this constant. At winding No. 5 the average winding-distance is already equal to $\pi$ to around 99,975 % !!

# Table 2 :

Analysis results referring to the distribution of the natural numbers ( or their square roots respectively ) on defined spiral-graph-systems on the "Square Root Spiral" ( --> see examples in FIG. 9-12 )

| **Number Group**<br>The Natural Numbers divisible by : | Number of Spiral Systems with a **positive** direction of rotation | Number of Spiral Systems with a **negative** direction of rotation | Naming of the Spiral Systems<br>( P = positive )<br>( N = negative ) | The " **2. Differential** " of the numbers on the Spiral Arms results in the following values :<br>( --> see examples on FIG.9, 10 | For every found spiral system **one** exemplary spiral arm is given :<br>**--> Specification of the number sequence** belonging to the chosen exemplary spiral arm<br>System | <br><br>Exemplary Spiral-arm ( -sequence ) |
|---|---|---|---|---|---|---|
| **19** | 1 | 1 | **N1**<br>**P1** | **19** | **N1**<br>**P1** | 19,76,152,247,361,....<br>19,38,76,133,209,.... |
| **17** | 1 | 1 | **N1**<br>**P1** | **17** | **N1**<br>**P1** | 51,136,238,357,....<br>85,136,204,289,.... |
| **13** | 1 | 2 | **N1 - N2**<br>**P1** | for N1 - N2 : **26**<br>for P1 : **13** | **N1**<br>**N2**<br>**P1** | 39,104,195,312,....<br>13,65,143,247,.....<br>39,65,104,156,..... |
| **11** | 2 | 2 | **N1 - N2**<br>**P1 - P2** | **22** | **N1**<br>**N2**<br>**P1**<br>**P2** | 22,77,154,253,374,....<br>11,55,121,209,319,....<br>33,44,77,132,209,....<br>33,55,99,165,253..... |
| **7** | 3 | 3 | **N1 - N3**<br>**P1 - P3** | **21** | **N1**<br>**N2**<br>**N3**<br>**P1**<br>**P2**<br>**P3** | 14,49,105,182,280,....<br>7,49,112,196,301,....<br>21,70,140,231,343,....<br>28,49,91,154,238,....<br>21,35,70,126,203,...<br>21,28,56,105,175,.... |
| **5** | 4 | 4 | **N1 - N4**<br>**P1 - P4** | **20** | **N1**<br>**N2**<br>**N3**<br>**N4**<br>**P1**<br>**P2**<br>**P3**<br>**P4** | 25,80,155,250,365,...<br>5,45,105,185,285,....<br>45,110,195,300,....<br>15,65,135,225,335,....<br>15,40,85,150,235,....<br>10,20,50,100,170,...<br>10,25,60,115,190,...<br>10,30,70,130,210,... |
| **3** | 6 | 7 | **N1 - N7**<br>**P1 - P6** | for N1 - N7 : **21**<br>for P1 - P6 : **18** | **N1**<br>**N2**<br>**N3**<br>**N4**<br>**N5**<br>**N6**<br>**N7**<br>**P1**<br>**P2**<br>**P3**<br>**P4**<br>**P5**<br>**P6** | 24,69,135,222,330,...<br>27,75,144,234,345,...<br>9,39,90,162,255,...<br>12,45,99,174,270<br>15,51,108,186,285,...<br>3,21,60,120,201,....<br>3,24,66,129,213,...<br>21,48,93,156,237,....<br>12,42,90,156,240,...<br>9,24,57,108,177,...<br>3,21,57,111,183,...<br>6,27,66,123,198,...<br>6,30,72,132,210,... |
| **2** | 9 | 10 | **N1 - N10**<br>**P1 - P9** | for N1 - N10 : **20**<br>for P1 - P9 : **18** | **N1**<br>**N2**<br>**N3**<br>**N4**<br>**N5**<br>**N6**<br>**N7**<br>**N8**<br>**N9**<br>**N10**<br>**P1**<br>**P2**<br>**P3**<br>**P4**<br>**P5**<br>**P6**<br>**P7**<br>**P8**<br>**P9** | 2,26,70,134,218,...<br>4,30,76,142,228,...<br>6,34,82,150,238,...<br>8,38,88,158,248,...<br>6,38,90,162,254,...<br>10,44,98,172,266,...<br>10,46,102,178,274,...<br>12,50,108,186,284,...<br>16,56,116,196,296,...<br>20,62,124,206,308,...<br>8,38,86,152,236,...<br>2,16,48,98,166,....<br>6,22,56,108,178,...<br>4,22,58,112,184,....<br>2,22,60,116,190,...<br>6,28,68,126,202,...<br>6,30,72,132,210,...<br>6,14,40,84,146,226,...<br>12,40,86,150,232,... |

| **Number Group** | Number of Spiral Arms with a **positive** direction of rotation | Number of Spiral Arms with a **negative** direction of rotation | Naming of the Spiral Arms<br>(Q = quadratic) | The " **2. Differential** " of the numbers on the Spiral Arms Q1, Q2 and Q3 results in the following values :<br>( --> see FIG. 1 ) | **Spiral Arm**<br>( FIG ) | **Number sequence belonging to Spiral Arm** |
|---|---|---|---|---|---|---|
| **The Square Numbers**<br>( 1, 4 , 9, 16, 25.....) | **3** | **None** | **Q1** to **Q3** | **18** | see FIG. 1 | **Q1 :** 1,16,49,100,169,256,361,484,.....<br>**Q2 :** 4,25,64,121,196,289,400,529,.....<br>**Q3 :** 9,36,81,144,225,324,441,576,..... |

**Table 3-A:** Quadratic Polynomials of exemplary Spiral-Graphs of the Number-Group-Spiral-Systems shown in FIG. 9-12 and described in Table 2

| **Numbers divisible by** | **2. Differential of Spiral-Graphs** | **Spiral Graph System** | **Number Sequence of one exemplary Spiral Graph of this system** | **Quadratic Polynomial 1** ( calculated with the first 3 numbers of the given sequence ) | **Quadratic Polynomial 2** ( calculated with 3 numbers starting with the **2.** Number of the sequence ) | **Quadratic Polynomial 3** ( calculated with 3 numbers starting with the **3.** Number of the sequence ) | **Quadratic Polynomial 4** ( calculated with 3 numbers starting with the **4.** Number of the sequence ) |
|---|---|---|---|---|---|---|---|
| **19** | **19** | N1 | 19 , 76 , 152 , 247 , 361 , 494 ,...... | $f_1 (x) = 9,5 x^2 + 28,5 x - 19$ | $f_2 (x) = 9,5 x^2 + 47,5 x + 19$ | $f_3 (x) = 9,5 x^2 + 66,5 x + 76$ | $f_4 (x) = 9,5 x^2 + 85,5 x + 152$ |
| | | P1 | 19 , 38 , 76 , 133 , 209 , 304 ,...... | $f_1 (x) = 9,5 x^2 - 9,5 x + 19$ | $f_2 (x) = 9,5 x^2 + 9,5 x + 19$ | $f_3 (x) = 9,5 x^2 + 28,5 x + 38$ | $f_4 (x) = 9,5 x^2 + 47,5 x + 76$ |
| **17** | **17** | N1 | 51 , 136 , 238 , 357 , 493 , 646 ,...... | $f_1 (x) = 8,5 x^2 + 59,5 x - 17$ | $f_2 (x) = 8,5 x^2 + 76,5 x + 51$ | $f_3 (x) = 8,5 x^2 + 93,5 x + 136$ | $f_4 (x) = 8,5 x^2 + 111 x + 238$ |
| | | P1 | 34 , 51 , 85 , 136 , 204 , 289 ,...... | $f_1 (x) = 8,5 x^2 - 8,5 x + 34$ | $f_2 (x) = 8,5 x^2 + 8,5 x + 34$ | $f_3 (x) = 8,5 x^2 + 25,5 x + 51$ | $f_4 (x) = 8,5 x^2 + 42,5 x + 85$ |
| **13** | **26** : for N1 - N2<br>**13** : for P1 | N1 | 39 , 104 , 195 , 312 , 455 , 624 ,...... | $f_1 (x) = 13 x^2 + 26 x + 0$ | $f_2 (x) = 13 x^2 + 52,0 x + 39$ | $f_3 (x) = 13 x^2 + 78,0 x + 104$ | $f_4 (x) = 13 x^2 + 104 x + 195$ |
| | | N2 | 13 , 65 , 143 , 247 , 377 , 533 ,...... | $f_1 (x) = 13 x^2 + 13 x - 13$ | $f_2 (x) = 13 x^2 + 39,0 x + 13$ | $f_3 (x) = 13 x^2 + 65,0 x + 65$ | $f_4 (x) = 13 x^2 + 91 x + 143$ |
| | | P1 | 26 , 39 , 65 , 104 , 156 , 221 ,...... | $f_1 (x) = 6,5 x^2 - 6,5 x + 26$ | $f_2 (x) = 6,5 x^2 + 6,5 x + 26$ | $f_3 (x) = 6,5 x^2 + 19,5 x + 39$ | $f_4 (x) = 6,5 x^2 + 32,5 x + 65$ |
| **11** | **22** | N1 | 22 , 77 , 154 , 253 , 374 , 517 ,...... | $f_1 (x) = 11 x^2 + 22 x - 11$ | $f_2 (x) = 11 x^2 + 44 x + 22$ | $f_3 (x) = 11 x^2 + 66 x + 77$ | $f_4 (x) = 11 x^2 + 88 x + 154$ |
| | | N2 | 11 , 55 , 121 , 209 , 319 , 451 ,...... | $f_1 (x) = 11 x^2 + 11 x - 11$ | $f_2 (x) = 11 x^2 + 33 x + 11$ | $f_3 (x) = 11 x^2 + 55 x + 55$ | $f_4 (x) = 11 x^2 + 77 x + 121$ |
| | | P1 | 33 , 44 , 77 , 132 , 209 , 308 ,...... | $f_1 (x) = 11 x^2 - 22 x + 44$ | $f_2 (x) = 11 x^2 + 0 x + 33$ | $f_3 (x) = 11 x^2 + 22 x + 44$ | $f_4 (x) = 11 x^2 + 44 x + 77$ |
| | | P2 | 33 , 55 , 99 , 165 , 253 , 363 ,...... | $f_1 (x) = 11 x^2 - 11 x + 33$ | $f_2 (x) = 11 x^2 + 11 x + 33$ | $f_3 (x) = 11 x^2 + 33 x + 55$ | $f_4 (x) = 11 x^2 + 55 x + 99$ |
| **7** | **21** | N1 | 14 , 49 , 105 , 182 , 280 , 399 ,...... | $f_1 (x) = 10,5 x^2 + 3,5 x + 0$ | $f_2 (x) = 10,5 x^2 + 24,5 x + 14$ | $f_3 (x) = 10,5 x^2 + 45,5 x + 49$ | $f_4 (x) = 10,5 x^2 + 66,5 x + 105$ |
| | | N2 | 7 , 49 , 112 , 196 , 301 , 427 ,...... | $f_1 (x) = 10,5 x^2 + 10,5 x - 14$ | $f_2 (x) = 10,5 x^2 + 31,5 x + 7$ | $f_3 (x) = 10,5 x^2 + 52,5 x + 49$ | $f_4 (x) = 10,5 x^2 + 73,5 x + 112$ |
| | | N3 | 21 , 70 , 140 , 231 , 343 , 476 ,...... | $f_1 (x) = 10,5 x^2 + 17,5 x - 7$ | $f_2 (x) = 10,5 x^2 + 38,5 x + 21$ | $f_3 (x) = 10,5 x^2 + 59,5 x + 70$ | $f_4 (x) = 10,5 x^2 + 80,5 x + 140$ |
| | | P1 | 28 , 49 , 91 , 154 , 238 , 343 ,...... | $f_1 (x) = 10,5 x^2 - 10,5 x + 28$ | $f_2 (x) = 10,5 x^2 + 10,5 x + 28$ | $f_3 (x) = 10,5 x^2 + 31,5 x + 49$ | $f_4 (x) = 10,5 x^2 + 52,5 x + 91$ |
| | | P2 | 21 , 35 , 70 , 126 , 203 , 301 ,...... | $f_1 (x) = 10,5 x^2 - 17,5 x + 28$ | $f_2 (x) = 10,5 x^2 + 3,5 x + 21$ | $f_3 (x) = 10,5 x^2 + 24,5 x + 35$ | $f_4 (x) = 10,5 x^2 + 45,5 x + 70$ |
| | | P3 | 21 , 28 , 56 , 105 , 175 , 266 ,...... | $f_1 (x) = 10,5 x^2 - 24,5 x + 35$ | $f_2 (x) = 10,5 x^2 + 3,5 x + 21$ | $f_3 (x) = 10,5 x^2 + 17,5 x + 28$ | $f_4 (x) = 10,5 x^2 + 38,5 x + 56$ |
| **5** | **20** | N1 | 25 , 80 , 155 , 250 , 365 , 500 ,...... | $f_1 (x) = 10 x^2 + 25 x - 10$ | $f_2 (x) = 10 x^2 + 45 x + 25$ | $f_3 (x) = 10 x^2 + 65 x + 80$ | $f_4 (x) = 10 x^2 + 85 x + 155$ |
| | | N2 | 5 , 45 , 105 , 185 , 285 , 405 ,...... | $f_1 (x) = 10 x^2 + 10 x - 15$ | $f_2 (x) = 10 x^2 + 30 x + 5$ | $f_3 (x) = 10 x^2 + 50 x + 45$ | $f_4 (x) = 10 x^2 + 70 x + 105$ |
| | | N3 | 45 , 110 , 195 , 300 , 425 , 570 ,...... | $f_1 (x) = 10 x^2 + 35 x + 0$ | $f_2 (x) = 10 x^2 + 55 x + 45$ | $f_3 (x) = 10 x^2 + 75 x + 110$ | $f_4 (x) = 10 x^2 + 95 x + 195$ |
| | | N4 | 15 , 65 , 135 , 225 , 335 , 465 ,...... | $f_1 (x) = 10 x^2 + 20 x + 15$ | $f_2 (x) = 10 x^2 + 40 x + 15$ | $f_3 (x) = 10 x^2 + 60 x + 65$ | $f_4 (x) = 10 x^2 + 80 x + 135$ |
| | | P1 | 15 , 40 , 85 , 150 , 235 , 340 ,...... | $f_1 (x) = 10 x^2 - 5 x + 10$ | $f_2 (x) = 10 x^2 + 15 x + 15$ | $f_3 (x) = 10 x^2 + 35 x + 40$ | $f_4 (x) = 10 x^2 + 55 x + 85$ |
| | | P2 | 10 , 20 , 50 , 100 , 170 , 260 ,...... | $f_1 (x) = 10 x^2 - 20 x + 20$ | $f_2 (x) = 10 x^2 + 0 x + 10$ | $f_3 (x) = 10 x^2 + 20 x + 20$ | $f_4 (x) = 10 x^2 + 40 x + 50$ |
| | | P3 | 10 , 25 , 60 , 115 , 190 , 285 ,...... | $f_1 (x) = 10 x^2 - 15 x + 15$ | $f_2 (x) = 10 x^2 + 5 x + 10$ | $f_3 (x) = 10 x^2 + 25 x + 25$ | $f_4 (x) = 10 x^2 + 45 x + 60$ |
| | | P4 | 10 , 30 , 70 , 130 , 210 , 310 ,...... | $f_1 (x) = 10 x^2 - 10 x + 10$ | $f_2 (x) = 10 x^2 + 10 x + 10$ | $f_3 (x) = 10 x^2 + 30 x + 30$ | $f_4 (x) = 10 x^2 + 50 x + 70$ |

**Table 3-B:** Quadratic Polynomials of exemplary Spiral-Graphs of the Number-Group-Spiral-Systems shown in FIG. 9-12 and described in Table 2

| Numbers divisible by | 2. Differential of Spiral-Graphs | Spiral Graph System | Number Sequence of <u>one</u> exemplary Spiral Graph of this system | Quadratic Polynomial 1 ( calculated with the first 3 numbers of the given sequence ) | Quadratic Polynomial 2 ( calculated with 3 numbers starting with the **2.** Number of the sequence ) | Quadratic Polynomial 3 ( calculated with 3 numbers starting with the **3.** Number of the sequence ) | Quadratic Polynomial 4 ( calculated with 3 numbers starting with the **4.** Number of the sequence ) |
|---|---|---|---|---|---|---|---|
| 3 | **21** : for N1 to N7<br>**18** : for P1 to P6 | N1 | 24 , 69 , 135 , 222 , 330 , 459 ,...... | $f_1 (x) = 10,5 x^2 + 13,5 x + 0$ | $f_2 (x) = 10,5 x^2 + 34,5 x + 24$ | $f_3 (x) = 10,5 x^2 + 55,5 x + 69$ | $f_4 (x) = 10,5 x^2 + 76,5 x + 135$ |
| | | N2 | 27 , 75 , 144 , 234 , 345 , 477 ,...... | $f_1 (x) = 10,5 x^2 + 16,5 x + 0$ | $f_2 (x) = 10,5 x^2 + 37,5 x + 27$ | $f_3 (x) = 10,5 x^2 + 58,5 x + 75$ | $f_4 (x) = 10,5 x^2 + 79,5 x + 144$ |
| | | N3 | 9 , 39 , 90 , 162 , 255 , 369 ,...... | $f_1 (x) = 10,5 x^2 - 1,5 x + 0$ | $f_2 (x) = 10,5 x^2 + 19,5 x + 9$ | $f_3 (x) = 10,5 x^2 + 40,5 x + 39$ | $f_4 (x) = 10,5 x^2 + 61,5 x + 90$ |
| | | N4 | 12 , 45 , 99 , 174 , 270 , 387 ,...... | $f_1 (x) = 10,5 x^2 + 1,5 x + 0$ | $f_2 (x) = 10,5 x^2 + 22,5 x + 12$ | $f_3 (x) = 10,5 x^2 + 43,5 x + 45$ | $f_4 (x) = 10,5 x^2 + 64,5 x + 99$ |
| | | N5 | 15 , 51 , 108 , 186 , 285 , 405 ,...... | $f_1 (x) = 10,5 x^2 + 4,5 x + 0$ | $f_2 (x) = 10,5 x^2 + 25,5 x + 15$ | $f_3 (x) = 10,5 x^2 + 46,5 x + 51$ | $f_4 (x) = 10,5 x^2 + 67,5 x + 108$ |
| | | N6 | 3 , 21 , 60 , 120 , 201 , 303 ,...... | $f_1 (x) = 10,5 x^2 - 13,5 x + 6$ | $f_2 (x) = 10,5 x^2 + 7,5 x + 3$ | $f_3 (x) = 10,5 x^2 + 28,5 x + 21$ | $f_4 (x) = 10,5 x^2 + 49,5 x + 60$ |
| | | N7 | 3 , 24 , 66 , 129 , 213 , 318 ,...... | $f_1 (x) = 10,5 x^2 - 10,5 x + 3$ | $f_2 (x) = 10,5 x^2 + 10,5 x + 3$ | $f_3 (x) = 10,5 x^2 + 31,5 x + 24$ | $f_4 (x) = 10,5 x^2 + 52,5 x + 66$ |
| | | P1 | 12 , 21 , 48 , 93 , 156 , 237 ,...... | $f_1 (x) = 9 x^2 - 18 x + 21$ | $f_2 (x) = 9 x^2 + 0 x + 12$ | $f_3 (x) = 9 x^2 + 18 x + 21$ | $f_4 (x) = 9 x^2 + 36 x + 48$ |
| | | P2 | 12 , 42 , 90 , 156 , 240 , 342 ,...... | $f_1 (x) = 9 x^2 + 3 x + 0$ | $f_2 (x) = 9 x^2 + 21 x + 12$ | $f_3 (x) = 9 x^2 + 39 x + 42$ | $f_4 (x) = 9 x^2 + 57 x + 90$ |
| | | P3 | 9 , 24 , 57 , 108 , 177 , 264 ,...... | $f_1 (x) = 9 x^2 - 12 x + 12$ | $f_2 (x) = 9 x^2 + 6 x + 9$ | $f_3 (x) = 9 x^2 + 24 x + 24$ | $f_4 (x) = 9 x^2 + 42 x + 57$ |
| | | P4 | 3 , 21 , 57 , 111 , 183 , 273 ,...... | $f_1 (x) = 9 x^2 - 9 x + 3$ | $f_2 (x) = 9 x^2 + 9 x + 3$ | $f_3 (x) = 9 x^2 + 27 x + 21$ | $f_4 (x) = 9 x^2 + 45 x + 57$ |
| | | P5 | 6 , 27 , 66 , 123 , 198 , 291 ,...... | $f_1 (x) = 9 x^2 - 6 x + 3$ | $f_2 (x) = 9 x^2 + 12 x + 6$ | $f_3 (x) = 9 x^2 + 30 x + 27$ | $f_4 (x) = 9 x^2 + 48 x + 66$ |
| | | P6 | 6 , 30 , 72 , 132 , 210 , 306 ,...... | $f_1 (x) = 9 x^2 - 3 x + 0$ | $f_2 (x) = 9 x^2 + 15 x + 6$ | $f_3 (x) = 9 x^2 + 33 x + 30$ | $f_4 (x) = 9 x^2 + 51 x + 72$ |
| 2 | **20** : for N1 to N10<br>**18** : for P1 to P9 | N1 | 2 , 26 , 70 , 134 , 218 , 322 ,...... | $f_1 (x) = 10 x^2 - 6 x - 2$ | $f_2 (x) = 10 x^2 + 14 x + 2$ | $f_3 (x) = 10 x^2 + 34 x + 26$ | $f_4 (x) = 10 x^2 + 54 x + 70$ |
| | | N2 | 4 , 30 , 76 , 142 , 228 , 334 ,...... | $f_1 (x) = 10 x^2 - 4 x - 2$ | $f_2 (x) = 10 x^2 + 16 x + 4$ | $f_3 (x) = 10 x^2 + 36 x + 30$ | $f_4 (x) = 10 x^2 + 56 x + 76$ |
| | | N3 | 6 , 34 , 82 , 150 , 238 , 346 ,...... | $f_1 (x) = 10 x^2 - 2 x - 2$ | $f_2 (x) = 10 x^2 + 18 x + 6$ | $f_3 (x) = 10 x^2 + 38 x + 34$ | $f_4 (x) = 10 x^2 + 58 x + 82$ |
| | | N4 | 8 , 38 , 88 , 158 , 248 , 358 ,...... | $f_1 (x) = 10 x^2 + 0 x - 2$ | $f_2 (x) = 10 x^2 + 20 x + 8$ | $f_3 (x) = 10 x^2 + 40 x + 38$ | $f_4 (x) = 10 x^2 + 60 x + 88$ |
| | | N5 | 6 , 38 , 90 , 162 , 254 , 366 ,...... | $f_1 (x) = 10 x^2 + 2 x - 6$ | $f_2 (x) = 10 x^2 + 22 x + 6$ | $f_3 (x) = 10 x^2 + 42 x + 38$ | $f_4 (x) = 10 x^2 + 62 x + 90$ |
| | | N6 | 10 , 44 , 98 , 172 , 266 , 380 ,...... | $f_1 (x) = 10 x^2 + 4 x - 4$ | $f_2 (x) = 10 x^2 + 24 x + 10$ | $f_3 (x) = 10 x^2 + 44 x + 44$ | $f_4 (x) = 10 x^2 + 64 x + 98$ |
| | | N7 | 10 , 46 , 102 , 178 , 274 , 390 ,...... | $f_1 (x) = 10 x^2 + 6 x - 6$ | $f_2 (x) = 10 x^2 + 26 x + 10$ | $f_3 (x) = 10 x^2 + 46 x + 46$ | $f_4 (x) = 10 x^2 + 66 x + 102$ |
| | | N8 | 12 , 50 , 108 , 186 , 284 , 402 ,...... | $f_1 (x) = 10 x^2 + 8 x - 6$ | $f_2 (x) = 10 x^2 + 28 x + 12$ | $f_3 (x) = 10 x^2 + 48 x + 50$ | $f_4 (x) = 10 x^2 + 68 x + 108$ |
| | | N9 | 16 , 56 , 116 , 196 , 296 , 416 ,...... | $f_1 (x) = 10 x^2 + 10 x - 4$ | $f_2 (x) = 10 x^2 + 30 x + 16$ | $f_3 (x) = 10 x^2 + 50 x + 56$ | $f_4 (x) = 10 x^2 + 70 x + 116$ |
| | | N10 | 20 , 62 , 124 , 206 , 308 , 430 ,...... | $f_1 (x) = 10 x^2 + 12 x - 2$ | $f_2 (x) = 10 x^2 + 32 x + 20$ | $f_3 (x) = 10 x^2 + 52 x + 62$ | $f_4 (x) = 10 x^2 + 72 x + 124$ |
| | | P1 | 8 , 38 , 86 , 152 , 236 , 338 ,...... | $f_1 (x) = 9 x^2 + 3 x - 4$ | $f_2 (x) = 9 x^2 + 21 x + 8$ | $f_3 (x) = 9 x^2 + 39 x + 38$ | $f_4 (x) = 9 x^2 + 57 x + 86$ |
| | | P2 | 2 , 16 , 48 , 98 , 166 , 252 ,...... | $f_1 (x) = 9 x^2 - 13 x + 6$ | $f_2 (x) = 9 x^2 + 5 x + 2$ | $f_3 (x) = 9 x^2 + 23 x + 16$ | $f_4 (x) = 9 x^2 + 41 x + 48$ |
| | | P3 | 6 , 22 , 56 , 108 , 178 , 266 ,...... | $f_1 (x) = 9 x^2 - 11 x + 8$ | $f_2 (x) = 9 x^2 + 7 x + 6$ | $f_3 (x) = 9 x^2 + 25 x + 22$ | $f_4 (x) = 9 x^2 + 43 x + 56$ |
| | | P4 | 4 , 22 , 58 , 112 , 184 , 274 ,...... | $f_1 (x) = 9 x^2 - 9 x + 4$ | $f_2 (x) = 9 x^2 + 9 x + 4$ | $f_3 (x) = 9 x^2 + 27 x + 22$ | $f_4 (x) = 9 x^2 + 45 x + 58$ |
| | | P5 | 2 , 22 , 60 , 116 , 190 , 282 ,...... | $f_1 (x) = 9 x^2 - 7 x + 0$ | $f_2 (x) = 9 x^2 + 11 x + 2$ | $f_3 (x) = 9 x^2 + 29 x + 22$ | $f_4 (x) = 9 x^2 + 47 x + 60$ |
| | | P6 | 6 , 28 , 68 , 126 , 202 , 296 ,...... | $f_1 (x) = 9 x^2 - 5 x + 2$ | $f_2 (x) = 9 x^2 + 13 x + 6$ | $f_3 (x) = 9 x^2 + 31 x + 28$ | $f_4 (x) = 9 x^2 + 49 x + 68$ |
| | | P7 | 6 , 30 , 72 , 132 , 210 , 306 ,...... | $f_1 (x) = 9 x^2 - 3 x + 0$ | $f_2 (x) = 9 x^2 + 15 x + 6$ | $f_3 (x) = 9 x^2 + 33 x + 30$ | $f_4 (x) = 9 x^2 + 51 x + 72$ |
| | | P8 | 6 , 14 , 40 , 84 , 146 , 226 ,...... | $f_1 (x) = 9 x^2 - 19 x + 16$ | $f_2 (x) = 9 x^2 - 1 x + 6$ | $f_3 (x) = 9 x^2 + 17 x + 14$ | $f_4 (x) = 9 x^2 + 35 x + 40$ |
| | | P9 | 12 , 40 , 86 , 150 , 232 , 332 ,...... | $f_1 (x) = 9 x^2 + 1 x + 2$ | $f_2 (x) = 9 x^2 + 19 x + 12$ | $f_3 (x) = 9 x^2 + 37 x + 40$ | $f_4 (x) = 9 x^2 + 55 x + 86$ |

| SQUARE NUMBERS | 2. Differential of Spiral-Graphs | SPIRAL GRAPH | Number Sequences of the Spiral Graphs Q1, Q2 and Q3 ( containing the Square Numbers ) | Quadratic Polynomial 1 ( calculated with the first 3 numbers of the given sequence ) | Quadratic Polynomial 2 ( calculated with 3 numbers starting with the **2.** Number of the sequence ) | Quadratic Polynomial 3 ( calculated with 3 numbers starting with the **3.** Number of the sequence ) | Quadratic Polynomial 4 ( calculated with 3 numbers starting with the **4.** Number of the sequence ) |
|---|---|---|---|---|---|---|---|
| **4, 9, 16, 25,..** | **18** | **Q1** | 1 , 16 , 49 , 100 , 169 , 256 ,...... | $f_1 (x) = 9 x^2 - 12 x + 4$ | $f_2 (x) = 9 x^2 + 6 x + 1$ | $f_3 (x) = 9 x^2 + 24 x + 16$ | $f_4 (x) = 9 x^2 + 42 x + 49$ |
| | | **Q2** | 4 , 25 , 64 , 121 , 196 , 289 ,...... | $f_1 (x) = 9 x^2 - 6 x + 1$ | $f_2 (x) = 9 x^2 + 12 x + 4$ | $f_3 (x) = 9 x^2 + 30 x + 25$ | $f_4 (x) = 9 x^2 + 48 x + 64$ |
| | | **Q3** | 9 , 36 , 81 , 144 , 225 , 324 ,...... | $f_1 (x) = 9 x^2 + 0 x + 0$ | $f_2 (x) = 9 x^2 + 18 x + 9$ | $f_3 (x) = 9 x^2 + 36 x + 36$ | $f_4 (x) = 9 x^2 + 54 x + 81$ |

**FIG. 17** : bare Square Root Spiral

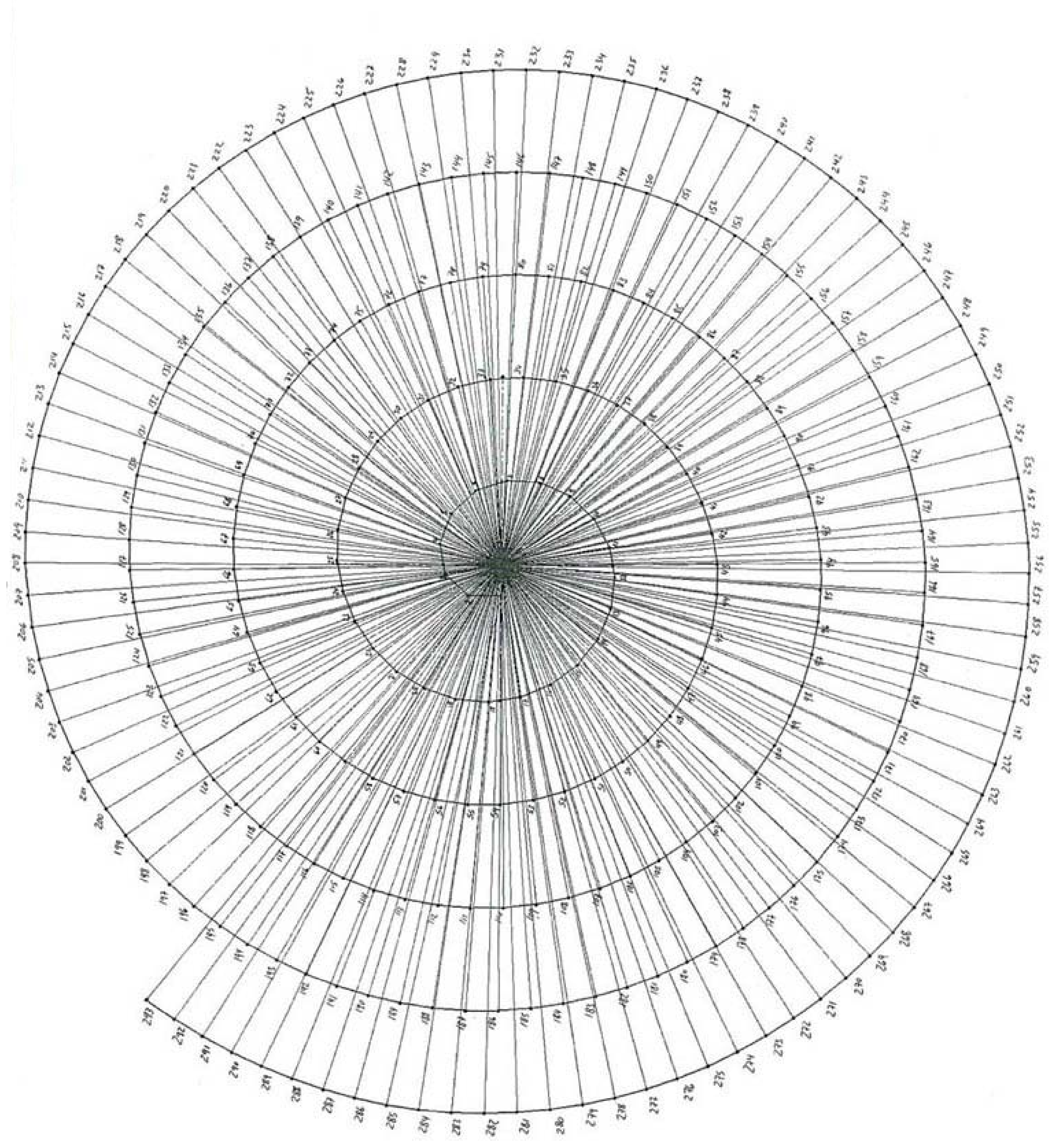